\documentclass[12pt,oneside,reqno]{amsart}
\usepackage{graphicx}
\usepackage{mathrsfs}
\usepackage{stmaryrd}
\usepackage{amsfonts}
\usepackage{enumerate,amsmath,amssymb,amsthm}
\newcommand{\be}{\begin{eqnarray}}
\newcommand{\ee}{\end{eqnarray}}
\newcommand{\ce}{\begin{eqnarray*}}
\newcommand{\de}{\end{eqnarray*}}
\newtheorem{theorem}{Theorem}[section]
\newtheorem{lemma}[theorem]{Lemma}
\newtheorem{remark}[theorem]{Remark}
\newtheorem{definition}[theorem]{Definition}
\newtheorem{proposition}[theorem]{Proposition}

\newtheorem{corollary}[theorem]{Corollary}
\def\e{\varepsilon}

\def\a{\alpha}
\def\om{\omega}
\def\Om{\Omega}
\def\b{\beta}
\def\p{\partial}

\def\l{\lambda}

\def\[{{\Big[}}
\def\]{{\Big]}}
\def\<{{\langle}}
\def\>{{\rangle}}
\def\({{\Big(}}
\def\){{\Big)}}

\def\dif{{\mathord{{\rm d}}}}

\def\no{\nonumber}
\def\bt{\begin{theorem}}
\def\et{\end{theorem}}
\def\bl{\begin{lemma}}
\def\el{\end{lemma}}
\def\br{\begin{remark}}
\def\er{\end{remark}}
\def\bd{\begin{definition}}
\def\ed{\end{definition}}
\def\bp{\begin{proposition}}
\def\ep{\end{proposition}}
\def\bc{\begin{corollary}}
\def\ec{\end{corollary}}

\def\cB{{\mathcal B}}

\def\cF{{\mathcal F}}

\def\cM{{\mathcal M}}

\def\cO{{\mathcal O}}

\def\cS{{\mathcal S}}

\def\mA{{\mathbb A}}

\def\mE{{\mathbb E}}
\def\mF{{\mathbb F}}

\def\mH{{\mathbb H}}

\def\mK{{\mathbb K}}

\def\mN{{\mathbb N}}

\def\mR{{\mathbb R}}

\def\mU{{\mathbb U}}

\def\mX{{\mathbb X}}
\def\mY{{\mathbb Y}}

\def\sD{{\mathscr D}}

\def\sH{{\mathscr H}}

\def\fM{{\mathfrak M}}
\def\fA{{\mathfrak A}}

\def\geq{\geqslant}
\def\leq{\leqslant}

\allowdisplaybreaks

\begin{document}
\title{On Stochastic Evolution Equations with non-Lipschitz Coefficients}
\date{}
\author{ Xicheng Zhang }

\dedicatory{
Department of Mathematics,
Huazhong University of Science and Technology,\\
Wuhan, Hubei 430074, P.R.China\\
Fakult\"at f\"ur Mathematik,
Universit\"at Bielefeld\\
Postfach 100131,
D-33501 Bielefeld, Germany\\
Email: XichengZhang@gmail.com
 }

\begin{abstract}

In this paper, we study the existence and uniqueness of solutions for
several classes of stochastic evolution equations with non-Lipschitz coefficients,
that is, backward stochastic evolution equations, stochastic Volterra type evolution equations
and stochastic functional evolution equations. In particular, the results can be used to treat a
large class of quasi-linear stochastic equations, which includes the reaction diffusion
and porous medium equations.
\end{abstract}



\maketitle
\allowdisplaybreaks

 \tableofcontents
\section{Introduction}

Let $\cO$ be a bounded open subset of $\mR^d$.
Consider the following stochastic porous medium equation
with Dirichlet boundary condition:
\be
\label{porous}\left\{
\begin{array}{ll}
\dif u_t=|w_t|\cdot\Delta(|u_t|^{p-2}u_t)\dif t+\dif w_t,\\
u_t(x)=0,\quad x\in\p\cO, t>0\\
u_0=\phi\in L^p(\cO),
\end{array}
\right.
\ee
where $p\geq 2$, $\Delta$ is the usual Laplace operator, and
$\{w_t,t\geq 0\}$ is a one dimensional standard Brownian motion. This is a degenerate
non-linear stochastic partial differential equation. Notice that the degeneracy may be caused by
$w_t=0$ and $u_t=0$. In the deterministic case, it is well known that
porous medium equations can be written as abstract monotone operator equations(cf. \cite{Ze} \cite{Sh}).
Thus, in the stochastic case, it can fall into a class of stochastic evolution equations studied
by Krylov-Rozovskii \cite{Kr-Ro}. More discussions about the stochastic porous medium equation
are referred to \cite{DR1} \cite{Re-Ro-Wa} \cite{Roe}.

On the other hand,  let us consider the following stochastic reaction diffusion equation:
\be
\label{rea}\left\{
\begin{array}{ll}
\dif u_t=|w_t|\cdot(\Delta u_t-|u_t|^{p-2}\cdot u_t)\dif t+\sqrt{|w_t|}\cdot u_t\dif w_t,\\
u_t(x)=0,\quad x\in\p\cO, t>0\\
u_0=\phi\in L^2(\cO),
\end{array}
\right.
\ee
where $p\geq 2$.
Usually, one wants to find an adapted process $u$ such that for almost all $w$
$$
u_\cdot(w)\in L^2([0,T],H^{1}_0(\cO))\cap L^p([0,T]\times\cO)\cap C([0,T], L^2(\cO)),
$$
and (\ref{rea}) holds in the generalized sense, where $H^1_0(\cO)$ is the usual Sobolev space.

However, from the well known results, it seems that one
cannot solve Eq.(\ref{porous}) and Eq.(\ref{rea}) because of the presence of
$|w_t|$ in front of the Laplace operator.
One of the main purposes in this paper is to extend the well known results in
\cite{Kr-Ro} \cite{Go-Mi} so that we can solve Eq.(\ref{porous}) and Eq.(\ref{rea}) in the generalized sense
for almost all path $w_t$.

In the present paper, we shall work on the framework of evolution triple. This is crucial for
treating a wide class of quasi-linear stochastic partial differential equations(including
reaction diffusion equations and porous medium equations).
We now recall some well-known results in this direction. In \cite{Pa1} \cite{Pa2}, Pardoux considered
linear stochastic partial differential equations(SPDEs) using the monotonicity method.
In \cite{Kr-Ro}, basing on their established It\^o's formula,
Krylov and Rozovskii proved a more general result under some monotonicity or dissipative conditions.
This classic work was later extended in several aspects:  to stochastic evolution equations(SEEs) driven by
general (discontinuous) martingales in \cite{Go1}, to SEEs with coercivity constants depending on $t$
in \cite{Go-Mi}, to SEEs related to some Orlicz spaces in \cite{Re-Ro-Wa}.
All these works are based on Galerkin's approximation. It should be remarked
that the semigroup method is another main tool in the theory of
semi-linear SPDEs (cf. \cite{DaZa} \cite{b1} \cite{b2} \cite{Kr1}
\cite{Hu-Le} \cite{Zh1} \cite{Zh3}etc.).
In order to solve Eq.(\ref{porous}), we need to deal with  SEEs
with random coercivity coefficients. This is our first
goal, and will be done in Section 3 after some preliminaries of Section 2.
Here, some stopping times techniques will be used.

The second aim is to prove the existence and uniqueness of solutions to
backward stochastic evolution equations.
Since Pardoux and Peng in \cite{Pa-Pe} proved the existence and uniqueness of solutions to
nonlinear backward stochastic differential equations(BSDEs), the theory of BSDEs has already been developed
extensively. It is well known that BSDEs can be applied to the studies of stochastic controls,
mathematics finances, deterministic PDEs, etc.. Meanwhile, backward SPDEs have also been studied in \cite{Hu-Ma-Yo} \cite{Ok-Pr-Zha} etc..
In these works, the authors mainly concentrated on semilinear BSPDEs.
The second aim in this paper is to prove the existence and
uniqueness of solutions to BSEEs with non-Lipschitz coefficients
in the framework of evolution triple.
Thus, it can be used to deal with a large class of quasi linear BSPDE.
We remark that Mao in \cite{Ma} has already studied the BSDEs with non-Lipschitz coefficients,
and the authors in \cite{Br-Ca} also investigated the BSDEs with monotone and arbitrary growth
coefficients.
This is the content of Section 4.

The third aim is to study the stochastic functional integral evolution equations
with non-Lipschitz coefficients,
which in particular includes a class of stochastic Volterra type evolution equations. Stochastic
Volterra equations driven by Brownian motion
were first studied by Berger-Mizel \cite{Be-Mi}. Later, Protter \cite{Pr} proved
the existence and uniqueness of stochastic Volterra equations driven by general semimartingales.
Recently, Wang in \cite{Wa-Zh} studied the the existence
and uniqueness of stochastic Volterra equations
with singular kernels and non-Lipschitz coefficients.
About the stochastic functional differential equations, Mohammend's book \cite{Mo} is
one of the main references.
In \cite{Ta-Li-Tr}, using the evolution semigroup approach,
the authors studied the existence, uniqueness and asymptotic behavior of
mild solutions to stochastic semilinear functional differential equations in Hilbert spaces.
In our proof of Section 5, the main tool is the usual Picard iteration.
As above, the results in Section 5 can be also used to deal with a class of quasi linear stochastic functional
partial differential equations.

Lastly, in Section 6 we shall present two applications for our abstract results:
stochastic porous medium equations
and stochastic reaction diffusion equations. In particular, Eq.(\ref{porous}) and Eq.(\ref{rea})
will be two special cases. It is worthy to say that the two examples given in
Section 6 have stochastic non-linear second order terms. Moreover,
we may also consider the corresponding backward and functional
stochastic partial differential equations with a slight modification.

\section{Framework and Preliminaries}

In this section we present a general setting in which we can deal with a large class of non-linear stochastic
partial differential equations, and also recall the powerful It\^o formula and a nonlinear Gronwall type inequality
(Bihari's inequality) for treating  non-Lipschitz equations.

Let $\mX$ be a reflexive and separable Banach space, which is densely injected in a separable Hilbert space
$\mH$. Identifying $\mH$ with its dual we get
$$
\mX\subset\mH\simeq\mH^*\subset\mX^*,
$$
where the star `$^*$' denotes the dual spaces.

Assume that the norm in $\mX$ is given by
$$
\|x\|_\mX:=\|x\|_{1,\mX}+\|x\|_{2,\mX},\quad x\in\mX.
$$
Denote by $\mX_i$, $i=1,2$ the completions of $\mX$ with respect to the norms $\|\cdot\|_{i,\mX}
=:\|\cdot\|_{\mX_i}$. Then $\mX=\mX_1\cap\mX_2$.
Let us also assume that both spaces are reflexive and embedded in $\mH$. Thus, we get two triples:
$$
\mX_1\subset\mH\simeq\mH^*\subset\mX^*_1, \ \ \mX_2\subset\mH\simeq\mH^*\subset\mX^*_2.
$$
Noticing that $\mX^*_1$ and $\mX^*_2$ can be thought as subspaces of $\mX^*$, one may define
a  Banach space $\mY:=\mX^*_1+\mX^*_2\subset\mX^*$ as follows: $f\in\mY$ if and only if
$f=f_1+f_2$, $f_i\in\mX^*_i, i=1,2$ and the norm of $f$ is defined by
$$
\|f\|_{\mY}=\inf_{f=f_1+f_2}(\|f_1\|_{\mX^*_1}+\|f_2\|_{\mX^*_2}).
$$
In the following, the dual pairs of $(\mX,\mX^*)$ and $(\mX_i,\mX^*_i), i=1,2$ are denoted respectively by
$$
[\cdot,\cdot]_{\mX},\quad[\cdot,\cdot]_{\mX_i}, \ \ \ i=1,2.
$$
Then, for any $x\in\mX$ and $f=f_1+f_2\in\mY\subset\mX^*$,
$$
[x,f]_{\mX}=[x,f_1]_{\mX_1}+[x,f_2]_{\mX_2}.
$$
We remark that if $f\in\mH$ and $x\in\mX$, then
$$
[x,f]_{\mX}=[x,f]_{\mX_1}=[x,f]_{\mX_2}=\<x,f\>_\mH,
$$
where $\<\cdot,\cdot\>_\mH$ stands for the inner product in $\mH$.

Let $(\Om,\cF,(\cF_t)_{t\geq 0}, P)$ be a complete separable filtration probability space, and
$Q$ a nonnegative definite and symmetric bounded linear operator on another Hilbert space $\mU$.
A cylindrical $Q$-Wiener process $\{W(t),t\geq0\}$ defined on
$(\Om,\cF,P)$  is given and assumed to be adapted to $(\cF_t)_{t\geq 0}$(cf. \cite{DaZa}).
In the following we shall only consider the case of $Q\equiv I$ for simplicity.
Let $L_2(\mU,\mH)$ denote the Hilbert space consisting of all the Hilbert-Schmidt operators
from $\mU$ to $\mH$, where the norm is denoted by $\|\cdot\|_{L_2(\mU,\mH)}$, and the inner product by
$\<\cdot,\cdot,\>_{L_2(\mU,\mH)}$.

Fix $T>0$. Let $\cM$ be the total of progressively measurable subsets of
$[0,T]\times\Omega$.
The following It\^o formula is taken from Gy\"ongy-Krylov \cite{Go-Ky1}.
\bt\label{Ito}
Let $X_0$ be an $\cF_0$-measurable $\mH$-valued random variable. Let
$$
Y_i:[0,T]\times\Omega\to\mX^*_i\in \cM/\cB(\mX^*_i), \ \ i=1,2,
$$
and $M$ an $\mH$-valued continuous locally square integrable martingale
starting form zero.
Let $\l_1,\l_2$ be two $\cM/\cB(\mR)$-measurable real valued processes such that
for $(\dif t\times\dif P)$-almost all $(t,\om)$, $\l_1(t,\om),\l_2(t,\om)>0$. Assume that for some $q_1,q_2>1$
and for almost all $\om$,
$$
\lambda_i(\cdot,\om)\in L^1([0,T],\dif t),\quad Y_i(\cdot,\om)
\cdot\l^{-\frac{1}{q_i}}_i(\cdot,\om)\in L^{\frac{q_i}{q_i-1}}([0,T],\dif t;\mX^*),\ \ i=1,2.
$$
Define an $\mX^*$-valued process by
\ce
X(t):=X_0+\int^t_0Y_1(s)\dif s+\int^t_0Y_2(s)\dif s+M(t).
\de
If there exists a $(\dif t\times\dif P)$-version $\tilde X$ of $X$ such that for almost all $\om$,
$$
\tilde X(\cdot,\om)\cdot\l^{\frac{1}{q_i}}_i(\cdot,\om)\in L^{q_i}([0,T],\dif t;\mX_i),\ \ i=1,2,
$$
then for almost all $\om$,
\begin{enumerate}[(i)]
\item $[0,T]\ni t\mapsto X(t,\om)\in\mH$ is continuous;

\item for all $t\in[0,T]$
\ce
\|X(t,\om)\|^2_\mH&=&\|X_0(\om)\|^2_\mH+2\int^t_0[\tilde X(s,\om),(Y_1+Y_2)(s,\om)]_\mX\dif s\\
&&+2\int^t_0\<X(s), \dif M(s)\>_\mH(\om)+\<M\>(t,\om),
\de
where $\<\cdot\>$ denotes the quadratic variation of $\mH$-valued  martingale.
\end{enumerate}
\et
\begin{proof}
Set $N_i(t):=\int^t_0\l^{\frac{1}{q_i}}_i(s)\dif s$ and
$\tilde Y_i(t):=Y_i(t)\cdot\l^{-\frac{1}{q_i}}_i(s), i=1,2$. Then
\ce
X(t)=X_0+\int^t_0\tilde Y_1(s)\dif N_1(s)+\int^t_0\tilde Y_2(s)\dif N_2(s)+M(t).
\de
By the assumptions and H\"older's inequality, we have for almost all $\om$,
\ce
\tilde Y_i(\cdot,\om)&\in& L^1([0,T],\dif N_i;\mX^*_i),\ \ i=1,2,\\
\tilde X(\cdot,\om)&\in& \cap_{i=1,2}L^1([0,T],\dif N_i;\mX_i).
\de
Moreover, by H\"older's inequality we have for $i=1,2$ and almost all $\omega$
\ce
&&\int^T_0\|\tilde Y_i(t,\om)\|_{\mX^*_i}\cdot\|\tilde X(t,\om)\|_{\mX_i}\dif N_i(t)\\
&=&\int^T_0\|Y_i(t,\om)\|_{\mX^*_i}\l^{-\frac{1}{q_i}}_i(t)\cdot
\|\tilde X(t,\om)\|_{\mX_i}\l^{\frac{1}{q_i}}_i(t)\dif t\\
&\leq&\left(\int^T_0\|Y_i(t,\om)\|^{\frac{q_i}{q_i-1}}_{\mX^*_i}\l^{-\frac{1}{q_i-1}}_i(t)
\dif t\right)^{\frac{q_i-1}{q_i}}\\
&&\times\left(\int^T_0\|\tilde X(t,\om)\|^{q_i}_{\mX_i}\l_i(t)\dif t\right)^{\frac{1}{q_i}}<+\infty.
\de
Thus, we can prove this Theorem along the same lines as in the proof of \cite[Theorem 2]{Go-Ky1}
(see also \cite{Kr-Ro} \cite{Ro} \cite{Roe}).
We omit the details.
\end{proof}
We now recall the following Bihari inequality(cf. \cite{Bi}). A multi-parameter version with jump
was proved in \cite{Zh-Zh}.
\bl\label{le1}
Let $\rho:\mR^+\to\mR^+$
be a continuous and non-decreasing function. Let $g(s)$ and $\l(s)$ be
two strictly positive functions on $\mR^+$ such that for some $g_0>0$
$$
g(t)\leq g_0+\int_0^t \l(s)\cdot\rho(g(s))\dif s,\quad t\geq0.
$$
If $\l$ is locally integrable, then
$$
g(t)\leq G^{-1}\left(G(g_0)+\int_0^t \l(s)\dif s\right),
$$
where $G(x):=\int_{x_0}^x \frac{1}{\rho(y)}\dif y$ is well defined for some
$x_0>0$, and $G^{-1}$ is the inverse function of $G$.

In particular, if $g_0=0$ and for some $\e>0$
\be
\int_{0}^\e\frac{1}{\rho(x)}\dif x=+\infty,\label{rho}
\ee
then $g(t)\equiv 0$.
\el
\br
The typical concave functions satisfying (\ref{rho}) are given by $\rho_k(x)$, $k=1,2,\cdots,$
\be
\label{FUN}
\rho_k(x):=\left\{
\begin{array}{lcl}
c_0 \cdot x\cdot \Pi_{j=1}^k\log^j x^{-1},&& x\leq\eta\\
c_0\cdot \eta\cdot \Pi_{j=1}^k\log^j \eta^{-1}+c_0\cdot\rho'_k(\eta-)\cdot (x-\eta), &&x>\eta,
\end{array}
\right.
\ee
where $\log^j x^{-1}:=\log\log\cdots\log x^{-1}$ and $c_0>0$, $0<\eta<1/e^k$.
\er

In the sequel, we use the following convention: $c_0, c_1, \cdots$ will
denote positive constants whose values may change in
different occasions. Moreover, the following Young  inequality will be used frequently:
Let $a,b>0$ and $\a,\b>1$ satisfying $\frac{1}{\a}+\frac{1}{\b}=1$, then for any $\e>0$
\be
ab\leq  \e a^\a +\frac{b^\b}{(\a\e)^{\b/\a} \b}.\label{Young}
\ee
For simplicity of notations, we also write
$$
\fA:=([0,T]\times\Omega,\cB([0,T])\times\cF,\dif t\times\dif P)
$$
and
$$
\fA_a:=([0,T]\times\Omega,\cM,\dif t\times\dif P).
$$

We now introduce three evolution  operators used in the present paper:
\ce
A_i: [0,T]\times\Om\times\mX_i\to\mX^*_i&\in& \cM\times\cB(\mX_i)/\cB(\mX^*_i),\ \ i=1,2,\\
B: [0,T]\times\Om\times\mX\to L_2(\mU,\mH)&\in&\cM\times\cB(\mX)/\cB(L_2(\mU,\mH)).
\de
In the following, for the sake of simplicity, we write
$$
A=A_1+A_2\in\mY\subset \mX^*.
$$
Assume that
\begin{enumerate}[(\bf $\mathbf{H}$1)]
\item (Hemicontinuity) For any $(t,\om)\in[0,T]\times\Om$ and $x,y,z\in\mX$, the mapping
$$
[0,1]\ni\e\mapsto [x,A(t,\om, y+\e z)]_\mX
$$
is continuous.

\item (Weak monotonicity) There exists $0\leq \lambda_0\in L^1(\fA)$ such that for all
$x,y\in\mX$ and $(t,\om)\in[0,T]\times\Om$
\ce
&&2[x-y,A(t,\om,x)-A(t,\om,y)]_\mX+\|B(t,\om,x)-B(t,\om,y)\|^2_{L_2(\mU,\mH)}\\
&&\qquad \qquad\qquad\qquad\qquad\qquad\ \leq \lambda_0(t,\om)\cdot\|x-y\|^2_\mH.
\de

\item (Weak coercivity) There exist $q_1,q_2\geq 2, c_1>0$ and positive functions
$\lambda_1,\lambda_2, \lambda_3,\xi\in L^1(\fA)$ satisfying that for almost all $(t,\om)$
\be
0\leq \lambda_0(t,\om)<c_1\cdot(\lambda_1(t,\om)\wedge\lambda_2(t,\om))\label{Con2}
\ee
and
\be
(t,\om)\mapsto\lambda_i(t,\om)\cdot e^{\frac{q_i-2}{2}
\int^t_0\lambda_0(s,\om)\dif s}\in L^1(\fA),\ \ i=1,2,\label{Con}
\ee
where $\lambda_0$ is same as in ($\mathbf{H2}$),
and such that for all $x\in\mX$ and $(t,\om)\in[0,T]\times\Om$
\ce
&&2[x,A(t,\om,x)]_\mX+\|B(t,\om,x)\|^2_{L_2(\mU,\mH)}\\
&\leq& -\sum_{i=1,2}\Big(\lambda_i(t,\om)\cdot\|x\|^{q_i}_{\mX_i}\Big)+\lambda_3(t,\om)\cdot\|x\|^2_\mH+\xi(t,\om).
\de

\item (Boundedness) There exist $c_{A_i}>0$ and $0\leq \eta_i\in L^{\frac{q_i}{q_i-1}}(\fA)$, $i=1,2$
such that for all $x\in\mX$ and $(t,\om)\in[0,T]\times\Om$
\ce
\|A_i(t,\om,x)\|_{\mX^*_i}\leq \eta_i(t,\om)\cdot\lambda^{\frac{1}{q_i}}_i(t,\om)
+c_{A_i}\cdot\lambda_i(t,\om)\cdot\|x\|^{q_i-1}_{\mX_i},\ \ i=1,2,
\de
where $q_1$ and $q_2$ are same as in ($\mathbf{H3}$).
\end{enumerate}

In order to emphasize $\lambda_i,\xi$ and $q_i,\eta_i$ below, we shall say that
$$
(A,B) \mbox{ satisfies } \sH(\lambda_0,\lambda_1,\lambda_2,\lambda_3,\xi,\eta_1,\eta_2, q_1,q_2).
$$
If there are no special declarations, we always suppose that $q_i\geq 2, c_1, c_{A_i}>0$,
$\eta_i\in L^{\frac{q_i}{q_i-1}}(\fA)$, $i=1,2$ and
$\lambda_i,\xi\in L^1(\fA),i=0,1,2,3$  are strictly positive functions, and
(\ref{Con2})-(\ref{Con}) hold.

\br\label{r1}

By $\mathbf{(H3)}$, $\mathbf{(H4)}$ and Young's inequality (\ref{Young}),
it follows that for any $x\in\mX$ and $(t,\om)\in[0,T]\times\Om$
\ce
\|B(t,\om,x)\|^2_{L_2(\mU,\mH)}
&\leq& 2\sum_{i=1,2}\Big(c_{A_i}\cdot\lambda_i(t,\om)\cdot\|x\|^{q_i}_{\mX_i}
+\eta_i(t,\om)\cdot\lambda^{\frac{1}{q_i}}_i(t,\om)\cdot\|x\|_{\mX_i}\Big)\\
&&+\lambda_3(t,\om)\cdot\|x\|^2_\mH+\xi(t,\om)\\
&\leq& \sum_{i=1,2}\Big(c_B\cdot\lambda_i(t,\om)\cdot\|x\|^{q_i}_{\mX_i}
+\eta_i^{\frac{q_i}{q_i-1}}(t,\om)\Big)\\
&&+\lambda_3(t,\om)\cdot\|x\|^2_\mH+\xi(t,\om),
\de
where $c_B>1$ only depends on $c_{A_i}$ and $q_i$, $i=1,2$.

\er

The following lemma is well known(cf. \cite{Kr-Ro}).
\bl\label{Le2}
Let $(A, 0)$ satisfy $\sH(0,\lambda_1,\lambda_2,\lambda_3,\xi,\eta_1,\eta_2, q_1,q_2)$, and
$0\leq \tau\leq T$  a bounded random variable. Let $X$ and $Y_i(i=1,2)$ be  respectively
$\mX$ and $\mX^*_i$-valued measurable processes with
$$
1_{[0,\tau]}(\cdot)\cdot\lambda^{\frac{1}{q_i-1}}_i\cdot X\in L^{q_i-1}(\fA;\mX_i),
\quad 1_{[0,\tau]}(\cdot)\cdot Y_i\in L^1(\fA;\mX^*_i),\ \ i=1,2.
$$
Assume that for any $\mX$-valued measurable process $\Phi$ satisfying
$$
1_{[0,\tau]}(\cdot)\cdot\lambda^{\frac{1}{q_i-1}}_i\cdot \Phi\in L^{q_i-1}(\fA;\mX_i), \ \ i=1,2,
$$
it holds
\be
\mE\left(\int^\tau_0[X(s)-\Phi(s),Y(s)-A(s,\Phi(s))]_\mX\dif s\right)\leq 0,\label{e2}
\ee
where $Y=Y_1+Y_2\in\mY\subset\mX^*$.

Then $Y(t,\omega)=A(t,\omega,X(t,\omega))$ for almost all
$(t,\omega)\in\{(t,\omega): t\in[0,\tau(\omega)]\}$.
\el
\begin{proof}
For any $\e\in(0,1)$ and $\mX$-valued bounded measurable process $\phi$, letting
$\Phi=X-\e \phi$ in (\ref{e2})
and dividing both sides by $\e$, we get
\ce
\mE\left(\int^{\tau}_0[\phi(s),Y(s)-A(s,X(s)-\e\phi(s))]_\mX\dif s\right)\leq 0.
\de
By $\mathbf{(H4)}$ and the assumptions, we have
\ce
1_{[0,\tau]}(\cdot)\cdot\Big(\|Y(\cdot)\|_{\mX^*}
+\sup_{\e\in(0,1)}\|A(\cdot,X(\cdot)-\e\phi(\cdot))\|_{\mX^*}\Big)\in L^1(\fA).
\de
Hence, by $\mathbf{(H1)}$ and the dominated convergence theorem
\ce
\mE\left(\int^{\tau}_0[\phi(s),Y(s)-A(s,X(s))]_\mX\dif s\right)\leq 0.
\de
By changing $\phi$ to $-\phi$ and the arbitrariness of $\phi$, we conclude that $Y=A(\cdot, X)$.
\end{proof}

The following lemma is simple and will be used in Section 4.
A short proof is provided here for the reader's convenience.
\bl\label{l4}
Let $(S,\cS)$ be a measurable space. Let $X:\mR^d\times S\to \mR^d$ be a measurable field.
Assume that for every $s\in S$, $\mR^d\ni x\mapsto X(x,s)\in\mR^d$ is a homeomorphism. Then, the
inverse $(x,s)\mapsto X^{-1}(x,s)$ is also a measurable field, i.e.: for each $x\in\mR^d$,
$X^{-1}(x,\cdot)$ is $\cS/\cB(\mR^d)$-measurable.
\el
\begin{proof}
Fix $x\in\mR^d$. It suffices to prove that for any bounded open  set $U\subset \mR^d$
\be
S_1:=\{s: X^{-1}(x,s)\in \bar U\}\in\cS,\label{es112}
\ee
where $\bar U$ denotes the closure of $U$ in $\mR^d$.

Let $Q$ be the set of rational points in $\mR^d$. Then
\be
S_1=\cap_{k=1}^\infty\cup_{y\in Q\cap  U}\{s: \|X(y,s)-x\|_{\mR^d}<1/k\}=:S_2.\label{es12}
\ee
In fact, if $s\in S_1$, then there is a $y\in\bar U$ such that $x=X(y,s)$. Since $U$
is open and $X(\cdot,s)$ is continuous, there exists a sequence $y_n\in U\cap Q$
such that $y_n\rightarrow y$ and $X(y_n,s)\rightarrow X(y,s)=x$. So, $s\in S_2$.
On the other hand, if $s\in S_2$, then there is a sequence $y_n\in U\cap Q$ such that
$\lim_{n\rightarrow\infty}\|X(y_n,s)-x\|_{\mR^d}=0$, and so $y_n\rightarrow X^{-1}(x,s)\in\bar U$.
(\ref{es112}) now follows from (\ref{es12}).
\end{proof}

\section{Stochastic Evolution Equations in Banach Spaces}

In this section, we consider the following stochastic evolution equation:
\be
\label{Eq1}\left\{
\begin{array}{ll}
\dif X(t)=A(t,X(t))\dif t+B(t,X(t))\dif W(t),\\
X(0)=X_0\in\mH,
\end{array}
\right.
\ee
where $(A, B)$ satisfies $\sH(\lambda_0,\lambda_1,\lambda_2,\lambda_3,\xi,\eta_1,\eta_2, q_1,q_2)$.
Here and after, one should keep in mind that
$A=A_1+A_2\in\mY\subset\mX^*$, where $A_1\in\mX^*_1, A_2\in\mX^*_2$.

Set
\be
H(t,\om):=\int^{t}_0\lambda_3(s,\om)\dif s,\label{He}
\ee
and define
\ce
\theta_t(\om):=\inf\{s\in[0,T]: H(s,\om)\geq t\}.
\de
Then $t\mapsto\theta_t$ are continuous stopping times,
 and $\theta_t\uparrow T$ as $t\uparrow\infty$. Here, $\inf\{\emptyset\}=T$ by convention.

Set for each $m\in\mN$
\ce
\mu^m(\dif t\times\dif\om):=1_{\{t\leq\theta_m(\om)\}}(\dif t\times\dif P),
\de
and define completed measurable spaces
$$
\fM^m:=\overline{([0,T]\times\Omega,\cB([0,T])\times\cF)}^{\mu^m}
$$
and
$$
\fM_a^m:=\overline{([0,T]\times\Omega,\cM)}^{\mu^m}.
$$
We introduce the following stochastic Banach spaces for later use:
for each $m\in\mN$
\ce
\mK^m_{1,i}&:=&L^{\frac{q_i}{q_i-1}}(\fM^m,\lambda_i^{-\frac{1}{q_i-1}}(t,\om)
\cdot\mu^m(\dif t\times\dif\om); \mX^*_i),\ \ i=1,2,\\
\mK^m_{2,i}&:=&L^{q_i}(\fM^m,\lambda_i(t,\om)\cdot\mu^m(\dif t\times\dif\om); \mX_i), \ \ i=1,2,\\
\mK^m_3&:=&L^2(\fM_a^m,\mu^m(\dif t\times\dif\om); L_2(\mU,\mH)),\\
\mK^m_4&:=&L^2(\fM^m,\mu^m(\dif t\times\dif\om); \mH),\\
\mK^m_5&:=&L^2(\fM^m,\lambda_3(t,\om)\cdot\mu^m(\dif t\times\dif\om); \mH),
\de
where the norms are defined in a natural manner, and
denoted  by $\|\cdot\|_{\mK}$, where $\mK$ stands for the above spaces. For instance,
\ce
\|X\|_{\mK^m_{2,i}}:=\left[\mE\left(\int^{\theta_m}_0\|X(t)\|^{q_i}_{\mX_i}
\cdot\lambda_i(t)\dif t\right)\right]^{1/q_i},\ \ i=1,2.
\de
\br\label{r3}
If $\lambda_3$ is non-random, then for some $m$ sufficiently large,
$\theta_m\equiv T.$
In this case, we shall omit the superscript `$m$' of $\mK^m$.
\er
We need the following lemma.
\bl\label{l1}
\begin{enumerate}[(i)]
\item $\mK^m_{i,j}, i,j=1,2$ and $\mK^m_3,\mK^m_4,\mK^m_5$ are separable and reflexive Banach spaces.

\item For any $Y\in\mK^m_{1,i}$, we have $\mE\left(\int^{\theta_m}_0\|Y(t)\|_{\mX^*_i}\dif t\right)
\leq c_0\cdot\|Y\|_{\mK^m_{1,i}}$, where $i=1$ or $2$.

\item Let $\{Y_n,n\in\mN\}$ weakly converge to $Y$ in $\mK^m_{1,i}$, then for any $X\in\mK^m_{2,i}$
\ce
\lim_{n\rightarrow\infty}\mE\left(\int^{\theta_m}_0[X(t), Y_n(t)]_{\mX_i}\dif t\right)
=\mE\left(\int^{\theta_m}_0[X(t), Y(t)]_{\mX_i}\dif t\right),
\de
where $i=1$ or $2$.

\item Let $\{X_n,n\in\mN\}$ weakly converge to $X$ in $\mK^m_{2,i}$, then for any $Y\in\mK^m_{1,i}$
\ce
\lim_{n\rightarrow\infty}\mE\left(\int^{\theta_m}_0[X_n(t),Y(t)]_{\mX_i}\dif t\right)
=\mE\left(\int^{\theta_m}_0[X(t),Y(t)]_{\mX_i}\dif t\right),
\de
where $i=1$ or $2$.
Moreover, if $\{X_n,n\in\mN\}$ also weakly converges to $\bar X$ in $\mK^m_5$,
then $\bar X(t,\om)=X(t,\om)$ for $\mu^m$-almost all
$(t,\om)$.

\item Define a linear operator from $\mK^m_3$ to $\mK^m_4$ as
\be
J(G):=\int^{\cdot\wedge\theta_m}_0G(s)\dif W(s),\label{JJ}
\ee
then $J$ is a continuous linear operator. In particular, $J$ is continuous
with respect to the weak topologies.
\end{enumerate}
\el
\begin{proof}
(i). It  follows from the separabilities and reflexivities of $\mX_i,\mX^*_i,i=1,2, $ and $\mH, L_2(\mU,\mH)$.

(ii). By H\"older's inequality we have
\ce
\mE\left(\int^{\theta_m}_0\|Y(t)\|_{\mX^*_i}\dif t\right)&=&
\mE\left(\int^{\theta_m}_0\|Y(t)\|_{\mX^*_i}\lambda^{-1/q_i}_2(t)\cdot\lambda^{1/q_i}_2(t)\dif t\right)\\
&\leq&\|Y\|_{\mK^m_{1,i}}\left(\int^T_0\mE(\lambda_i(t))\dif t\right)^{1/q_i}.
\de

(iii). It follows from
\ce
X(\cdot)\cdot\lambda^{\frac{1}{q_i-1}}_i(\cdot)\in
L^{q_i}(\fM^m,\lambda^{-\frac{1}{q_i-1}}_i(t,\om)\cdot\mu^m(\dif t\times\dif \om); \mX_i)\subset (\mK^m_{1,i})^*.
\de

(iv). The first conclusion follows from
\ce
Y(\cdot)\cdot\lambda^{-1}_i(\cdot)\in
L^{\frac{q_i}{q_i-1}}(\fM^m,\lambda_i(t,\om)\cdot\mu^m(\dif t\times\dif \om); \mX^*_i)\subset (\mK^m_{2,i})^*.
\de
As for the second conclusion, by the well known Banach-Saks-Kakutani theorem,
there exists a subsequence of $X_n$(still denoted by $X_n$) such that
its C\'esaro means $\tilde X_n$ strongly converges to $X$ and $\bar X$
in $\mK^m_{2,i}$ and $\mK^m_5$ respectively.
Therefore, there is a subsequence $\tilde X_{n_k}$ such that
for $\mu^m$-almost all $(t,\om)$, $\tilde X_{n_k}(t,\om)\rightarrow X(t,\om)$ in $\mX$, and
$\tilde X_{n_k}(t,\om)\rightarrow \bar X(t,\om)$ in $\mH$. Since $\mX$ is continuously and densely
embedded in $\mH$, we have
$\bar X(t,\om)=X(t,\om)$ for $\mu^m$-almost all $(t,\om)$.

(v). It follows from
\ce
\|J(G)\|^2_{\mK^m_4}&=&\mE\left(\int^{\theta_m}_0\left\|\int^{t\wedge\theta_m}_0G(s)\dif W(s)\right\|^2_\mH\dif t\right)\\
&\leq&\int^{T}_0\mE\left(\int^{t\wedge\theta_m}_0\|G(s)\|^2_{L_2(\mU,\mH)}\dif s\right)\dif t\\
&\leq&T\|G\|^2_{\mK^m_3}.
\de
The proof is complete.
\end{proof}
\bd\label{d2}
An $\mH$-valued continuous $\cF_t$-adapted process $X(t,\om)$
is called a solution of Eq.(\ref{Eq1}) if for almost all $\omega\in\Om$,
$$
X(\cdot,\om)\in \cap_{i=1,2}L^{q_i}([0,T],\lambda_i(\cdot,\om)\dif t; \mX_i)
$$
and for all $t\in[0,T]$
$$
X(t,\om)=X_0(\om)+\int^t_0A(s,\om,X(s,\om))\dif s+\int^t_0B(s, X(s))\dif W(s)(\om),
$$
where the first integral is understood as an $\mX^*$-valued Bochner integral.
\ed
\br
Note that
\ce
\int^t_0A(s,\om,X(s,\om))\dif s=\int^t_0A_1(s,\om,X(s,\om))\dif s+\int^t_0A_2(s,\om,X(s,\om))\dif s.
\de
Since $X$ is $\cM/\cB(\mH)$-measurable, $1_{\mX_i}(X)\cdot X$ is  $\cM/\cB(\mX_i)$-measurable by
\cite[Lemma 2.1]{Kr-Ro} for $i=1,2$. The above integrals are meaningful.
\er
We have the following estimates for the solutions of Eq.(\ref{Eq1}).
\bt\label{th2}
Assume that $\mathbf{(H1)}$-$\mathbf{(H4)}$ hold and $X_0\in L^2(\Omega,\cF_0,P;\mH)$.
Let $X(t)$ be any solution of Eq. (\ref{Eq1}) in the sense of Definition \ref{d2}.
Then, we have for any $m\in\mN$
\ce
\mE\|X(\theta_m)\|^2_\mH+\sum_{i=1,2}\|X\|^{q_i}_{\mK^m_{2,i}}+\|X\|^2_{\mK^m_5}
\leq c_m\left(\mE\|X_0\|^2_\mH+\int^T_0\mE(\xi(s))\dif s\right),
\de
and
\ce
&&\mE\left(\sup_{t\in[0,\theta_m]}\|X(t)\|^2_\mH\right)+\|X\|^2_{\mK^m_4}
+\|B(\cdot,X(\cdot))\|^2_{\mK^m_3}+\sum_{i=1,2}\|A_i(\cdot,X(\cdot))\|^{\frac{q_i}{q_i-1}}_{\mK^m_{1,i}}\\
&\leq& c_m\left(\mE\|X_0\|^2_\mH+\int^T_0\mE\left(\xi(s)+\eta^{\frac{q_1}{q_1-1}}_1(s)
+\eta^{\frac{q_2}{q_2-1}}_2(s)\right)\dif s\right),
\de
where $c_m$  only  depends on $m$, $T$ and $c_{A_i}, q_i, i=1,2$.
\et
\begin{proof}
By It\^o's formula (Theorem \ref{Ito}) and $\mathbf{(H3)}$, we have
\be
&&\|X(t)\|^2_\mH-\|X_0\|^2_\mH-M(t)\no\\
&=&\int^{t}_0\left(2[X(s),A(s,X(s))]_\mX+\|B(s,X(s))\|^2_{L_2(\mU,\mH)}\right)\dif s\no\\
&\leq&\int^{t}_0\left(-\sum_{i=1,2}\Big(\lambda_i(s)\cdot\|X(s)\|^{q_i}_{\mX_i}\Big)
+\lambda_3(s)\cdot\|X(s)\|^2_\mH+\xi(s)\right)\dif s,\label{e5}
\ee
where $M(t)$ is a continuous local martingale given by
\ce
M(t):=2\int^{t}_0\<X(s), B(s,X(s))\dif W(s)\>_\mH.
\de

For any $R>0$, define the stopping time
\be
\tau_R:=\inf\left\{t\in[0,T]: \|X(t)\|_\mH\geq R, \int^t_0\lambda_i(s)\cdot
\|X(s)\|^{q_i}_{\mX_i}\dif s\geq R, i=1,2\right\}.\label{sto}
\ee
Then, by Definition \ref{d2}, $\tau_R\uparrow T$ a.s. as $R\uparrow\infty$.

By Remark \ref{r1} and the change of clock(cf. \cite{Re-Yo}),
we know that $\{M(\theta_t\wedge\tau_R), t\geq 0\}$
is a continuous $\cF_{\theta_t}$-martingale.
Indeed, this follows from
\ce
\<M(\theta_\cdot\wedge\tau_R)\>(t)\leq 4\int^{\theta_t\wedge\tau_R}_0\|X(s)\|^2_\mH
\cdot\|B(s,X(s))\|^2_{L_2(\mU,\mH)}\dif s\leq c_R.
\de
So, replacing $t$ by $\theta_t\wedge\tau_R$ in (\ref{e5})
and taking expectations for both sides of (\ref{e5}), we obtain
\ce
&&\mE\|X(\theta_t\wedge\tau_R)\|^2_\mH-\mE\|X_0\|^2_\mH+\sum_{i=1,2}\mE\left(\int^{\theta_t\wedge\tau_R}_0
\lambda_i(s)\cdot\|X(s)\|^{q_i}_{\mX_i}\dif s\right)\\
&\leq& \mE\left(\int^{\theta_t\wedge\tau_R}_0
\left(\lambda_3(s)\cdot\|X(s)\|^2_\mH+\xi(s)\right)\dif s\right)\\
&=&\mE\left(\int^{\theta_t\wedge\tau_R}_0\|X(s)\|^2_\mH\dif H(s)\right)
+\mE\left(\int^{\theta_t\wedge\tau_R}_0\xi(s)\dif s\right)\\
&\leq&\mE\left(\int^{\theta_t}_0\|X(s\wedge\tau_R)\|^2_\mH\dif H(s)\right)
+\int^T_0\mE(\xi(s))\dif s\\
&=&\int^{t}_0\mE\|X(\theta_s\wedge\tau_R)\|^2_\mH\dif s+\int^T_0\mE(\xi(s))\dif s,
\de
where $H(s)$ is defined by (\ref{He}), and
in the last step we have used the variable substitution formula.

Hence, by Gronwall's inequality we have for any $t\geq 0$
\ce
\mE\|X(\theta_t\wedge\tau_R)\|^2_\mH\leq e^{t}\left(\mE\|X_0\|^2_\mH+\int^T_0\mE(\xi(s))\dif s\right).
\de
Letting $R\rightarrow\infty$, by Fatou's lemma  we obtain that for any $m\in\mN$
\ce
\mE\|X(\theta_m)\|^2_\mH\leq e^m\left(\mE\|X_0\|^2_\mH+\int^T_0\mE(\xi(s))\dif s\right)
\de
as well as
\be
&&\sum_{i=1,2}\mE\left(\int^{\theta_m}_0
\lambda_i(s)\cdot\|X(s)\|^{q_i}_{\mX_i}\dif s\right)
+\mE\left(\int^{\theta_m}_0\lambda_3(s)\cdot\|X(s)\|^2_\mH\dif s\right)\no\\
&\leq& c_m\left(\mE\|X_0\|^2_\mH+\int^T_0\mE(\xi(s))\dif s\right),\label{e14}
\ee
which gives the first estimate.

From (\ref{e5}), by Burkholder's inequality and Young's inequality (\ref{Young}) we further have
\ce
&&\mE\left(\sup_{t\in[0,\theta_m]}\|X(t)\|^2_\mH\right)-\mE\|X_0\|^2_\mH\\
&\leq& \mE\left(\int^{\theta_m}_0\left(\lambda_3(s)\cdot\|X(s)\|^2_\mH
+\xi(s)\right)\dif s\right)\\
&&+c_0\mE\left(\int^{\theta_m}_0\|X(s)\|^2_\mH\cdot\|B(s,X(s))\|^2_{L_2(\mU,\mH)}\dif s\right)^{1/2}\\
&\leq&\int^T_0\mE(\xi(s))\dif s+\mE\left(\int^{\theta_m}_0
\lambda_3(s)\cdot\|X(s)\|^2_\mH\dif s\right)\\
&&+\frac{1}{2}\mE\left(\sup_{t\in[0,\theta_m]}\|X(t)\|^2_\mH\right)
+c_0\mE\left(\int^{\theta_m}_0\|B(s,X(s))\|^2_{L_2(\mU,\mH)}\dif s\right).
\de
Hence
\ce
\mE\left(\sup_{t\in[0,\theta_m]}\|X(t)\|^2_\mH\right)
&\leq&c_m\left(\mE\|X_0\|^2_\mH+\int^T_0\mE(\xi(s))\dif s\right)\\
&&+c_0\mE\left(\int^{\theta_m}_0\|B(s,X(s))\|^2_{L_2(\mU,\mH)}\dif s\right).
\de
The second estimate now follows from $\mathbf{(H4)}$, Remark \ref{r1} and (\ref{e14}).
\end{proof}
We  now prove our main result in this section.
\bt\label{th1}
Assume that $(A, B)$ satisfies $\sH(\lambda_0,\lambda_1,\lambda_2,\lambda_3,\xi,\eta_1,\eta_2, q_1,q_2)$.
Then for any $X_0\in L^2(\Omega,\cF_0,P;\mH)$,
there exists a unique solution to Eq.(\ref{Eq1}) in the sense of Definition \ref{d2}.
\et
\begin{proof}

(Uniqueness)

Let $X_1$ and $X_2$ be two solutions of Eq.(\ref{Eq1}) in the sense of Definition \ref{d2}.
For $t\geq 0$, define
$$
\beta_t:=\inf\left\{s\in[0,T]: \int^{s}_0\lambda_0(r)\dif r\geq t\right\},
$$
and for $R>0$ and $i=1,2$, let $\tau^i_R$ be defined as in (\ref{sto}) corresponding to $X_i$.
For $t_0\in(0,T)$, set
$\tau^{t_0}_R:=\tau^1_R\wedge\tau^2_R\wedge t_0$.
By It\^o's formula(Theorem \ref{Ito}), as in the proof of Theorem \ref{th2} we have
\ce
&&\mE\|(X_1-X_2)(\beta_t\wedge \tau^{t_0}_R)\|^2_\mH\\
&=&\mE\Bigg(\int^{\beta_t\wedge \tau^{t_0}_R}_0\Big(2[X_1(s)-X_2(s),A(s,X_1(s))-A(s,X_2(s))]_\mX\\
&&+\|B(s,X_1(s))-B(s,X_2(s))\|^2_{L_2(\mU,\mH)}\Big)\dif s\Bigg)\\
&\leq&\mE\Bigg(\int^{\beta_t\wedge \tau^{t_0}_R}_0\|X_1(s)-X_2(s)\|^2_\mH
\cdot\lambda_0(s)\dif s\Bigg)\\
&\leq&\mE\Bigg(\int^{t}_0\|(X_1-X_2)(\beta_s\wedge \tau^{t_0}_R)\|^2_\mH\dif s\Bigg).
\de
Using Gronwall's inequality yields that for any $t\geq 0$ and $R>0$
\ce
\mE\|(X_1-X_2)(\beta_t\wedge \tau^{t_0}_R)\|^2_\mH=0.
\de
Letting $R,t\rightarrow\infty$, and by Fatou's lemma we get
\ce
\mE\|(X_1-X_2)(t_0)\|^2_\mH=0.
\de
The uniqueness is then obtained.

(Existence)
 We divide the proof into five steps.

{\it (Step 1)}

First of all, let us   reduce $\mathbf{(H1)}$-$\mathbf{(H4)}$ to the case of $\lambda_0=0$.
Let $X$ be a solution of Eq.(\ref{Eq1}). Set
\ce
\gamma(t,\om)&:=&e^{\frac{1}{2}\int^t_0\lambda_0(s,\om)\dif s}\\
\tilde X(t,\om)&:=&\gamma^{-1}(t,\om)\cdot X(t,\om)\\
\tilde A_i(t,\om,x)&:=&\gamma^{-1}(t,\om)\cdot A_i(t,\om,\gamma(t,\om)x)-\frac{\lambda_0(t,\om)x}{2},\ \ i=1,2,\\
\tilde B(t,\om,x)&:=&\gamma^{-1}(t,\om)\cdot B(t,\om,\gamma(t,\om)x).
\de
Noticing that $\gamma^{-1}(t)\leq 1$, we have
\ce
&&2[x,\tilde A(t,x)]_\mX+\|\tilde B(t,x)\|^2_{L_2(\mU,\mH)}\\
&\leq&\gamma^{-2}(t)
\left(-\sum_{i=1,2}\lambda_i(t)\cdot \gamma^{q_i}(t)\cdot \|x\|^{q_i}_{\mX_i}+\lambda_3(t)\cdot
\gamma^2(t)\cdot \|x\|^2_\mH+\xi(t)\right)+\lambda_0(t)\cdot\|x\|^2_\mH\\
&\leq&-\sum_{i=1,2}\Big(\lambda_i(t)\cdot \gamma^{q_i-2}(t)\cdot \|x\|^{q_i}_{\mX_i}\Big)
+(\lambda_3(t)+\lambda_0(t))\cdot\|x\|^2_\mH+\xi(t),
\de
On the other hand, by  (\ref{Con2}) and Young's inequality (\ref{Young}) and $q_1,q_2\geq 2$,
we have for $i=1,2$
\ce
\|\tilde A_i(t,x)\|_{\mX^*_i}&\leq&\gamma^{-1}(t)\cdot \|A_i(t,\gamma(t)x)\|_{\mX^*_i}
+\frac{\lambda_0(t)\cdot\|x\|_{\mX^*_i}}{2}\\
&\leq&\gamma^{-1}(t)\left(\eta_i(t)\cdot \lambda^{1/q_i}_i(t)
+c_{A_i}\gamma^{q_i-1}(t)\cdot \lambda_i(t)\cdot \|x\|^{q_i-1}_{\mX_i}\right)
+c_0\lambda_i(t)\cdot \|x\|_{\mX_i}\\
&\leq&\left(\eta_i(t)+c_0\lambda^{\frac{q_i-1}{q_i}}_i(t)\right)\lambda^{1/q_i}_i(t)
+\left(c_{A_i}\gamma^{q_i-2}(t)+c_0\right)\lambda_i(t)\cdot \|x\|^{q_i-1}_{\mX_i}\\
&\leq&\left(\eta_i(t)+c_0\lambda^{\frac{q_i-1}{q_i}}_i(t)\right)(\gamma^{q_i-2}(t)\cdot \lambda_i(t))^{1/q_i}
+c_0\gamma^{q_i-2}(t)\cdot \lambda_i(t)\cdot \|x\|^{q_i-1}_{\mX_i}.
\de
Thus,  $(\tilde A, \tilde B)$ satisfies
$\sH(0,\tilde\lambda_1,\tilde\lambda_2,\lambda_3+\lambda_0,\xi,\tilde\eta_1,\tilde\eta_2, q_1,q_2)$  with
\ce
\tilde\lambda_i(t)&:=&\lambda_i(t)\cdot e^{\frac{q_i-2}{2}\int^t_0\lambda_0(s)\dif s}\in L^1(\fA),\ \ i=1,2,\\
\tilde\eta_i(t)&:=&\eta_i(t)+c_0\lambda^{\frac{q_i-1}{q_i}}_i(t)\in L^{\frac{q_i}{q_i-1}}(\fA),\ \ i=1,2.
\de
Moreover, it holds that
\ce
\dif \tilde X(t)=\tilde A(t,\tilde X(t))\dif t+\tilde B(t,\tilde X(t))\dif W(t),\quad  \tilde X(0)=X(0).
\de
It is easy to see that
$$
X(\cdot,\om)\in \cap_{i=1,2}L^{q_i}([0,T],\lambda_i(\cdot,\om)\dif t; \mX_i)$$
is equivalent to
$$\tilde X(\cdot,\om)\in \cap_{i=1,2}L^{q_i}([0,T],\tilde\lambda_i(\cdot,\om)\dif t; \mX_i).$$
So, we may assume that $\lambda_0=0$ in the following proof.

{\it (Step 2)}

We now use Galerkin's approximation to prove the existence of  solutions.

Let  $\{e_i,i\in\mN\}\subset\mX$ be a normal orthogonal basis of $\mH$. Set
$$
\Pi_nx:=\sum_{i=1}^n[e_i,x]_\mX \cdot e_i,\quad x\in\mX^*.
$$
Then, the mapping $\Pi_n: \mX^*\to\mX$ is linear and continuous, and satisfy
$$
\Pi_nx=\sum_{i=1}^n\<e_i,x\>_\mH \cdot  e_i, \quad x\in\mH
$$
and
$$
[\Pi_nx, y]_\mX=[\Pi_n y, x]_\mX, \quad x,y\in\mX^*.
$$

We also fix a normal orthogonal basis $\{f_1,f_2,\cdots\}$ of $\mU$.
Let $W_j(t):=\<W(t),f_j\>_\mU$ for $j\in\mN$.
Consider the following It\^o type stochastic ordinary differential equation in $\mR^n$:
\be
\left\{
\begin{array}{ll}
\dif X^i_n(t)=b^i(t, X_n(t))\dif t+\sum_{j=1}^n\sigma^i_j(t, X_n(t)) \dif W_j(t),\\
X^i_n(0)=\<X_0,e_i\>_\mH, \quad i=1,\cdots, n,
\end{array}
\right.\label{es6}
\ee
where $b^i(t, x):=[e_i, A(t,x\cdot e)]_\mX$, and $\sigma^i_j(t,x):=\<e_i,B(t,x\cdot e)(f_j)\>_\mH$.
Here, $x\in\mR^n$ and $x\cdot e:=\sum_{i=1}^n x^i e_i$.

The coefficients satisfy the following conditions by $\mathbf{(H1)}$-$\mathbf{(H4)}$:

(i) $b$ and $\sigma$ are $\cM\times\cB(\mR^n)/\cB(\mR^n)$ and
$\cM\times\cB(\mR^n)/\cB(L_2(\mR^n,\mR^n))$-measurable respectively and continuous in $x$.

(ii) For any $(t,\om)\in [0,T]\times\Omega$ and $x,y\in\mR^n$
\ce
2\<x-y,b(t,\om,x)-b(t,\om,y)\>_{\mR^n}+\|\sigma(t,\om,x)-\sigma(t,\om,y)\|^2_{L_2(\mR^n,\mR^n)}\leq 0.
\de

(iii) For any $(t,\om)\in [0,T]\times\Omega$ and $x\in\mR^n$
\ce
2\<x,b(t,\om,x)\>_{\mR^n}+\|\sigma(t,\om,x)\|^2_{L_2(\mR^n,\mR^n)}\leq
\lambda_3(t,\om)\cdot\|x\|^2_{\mR^n}+\xi(t,\om).
\de

(iv) For any $(t,\om)\in [0,T]\times\Omega$ and $x\in\mR^n$
\ce
\|b(t,\om,x)\|_{\mR^n}\leq
c_n\sum_{i=1,2}\left(\eta_i(t,\om)\cdot\lambda^{1/q_i}_i(t,\om)+\lambda_i(t,\om)\cdot\|x\|^{q_i-1}_{\mR^n}\right).
\de
By the well-known result (cf. \cite{Kr0}), there exists a
unique continuous $\cF_t$-adapted solution denoted by $X^i_n(t)$ to Eq.(\ref{es6}). Moreover,
if we let $X_n(t):=\sum_{i=1}^nX^i_n(t)e_i$, then we can write  Eq.(\ref{es6}) as
\be
X_n(t)=\Pi_nX_0+\int^t_0\Pi_nA(s,X_n(s))\dif s+\int^t_0\Pi_nB(s,X_n(s))\tilde\Pi_n\dif W(s),\label{e6}
\ee
where $\tilde\Pi_n$ is the projection on $\mbox{span}\{f_1,\cdots,f_n\}$ in $\mU$.

Noticing that
$$
\|\Pi_nB(s,X_n(s))\tilde\Pi_n\|_{L_2(\mU,\mH)}\leq\|B(s,X_n(s))\|_{L_2(\mU,\mH)},
$$
and using the same method as in the proof of Theorem \ref{th2},
by $\mathbf{(H4)}$ and Remark \ref{r1}, we have for all $n\in\mN$
\ce
&&\mE\|X_n(\theta_m)\|^2_\mH+\|X_n\|^2_{\mK^m_4}+\|X_n\|^2_{\mK^m_5}
+\|B(\cdot,X_n)\|^2_{\mK^m_3}\\
&&+\sum_{i=1,2}\Big(\|X_n\|^{q_i}_{\mK^m_{2,i}}+\|A_i(\cdot,X_n)\|^{\frac{q_i}{q_i-1}}_{\mK^m_i}\Big)\\
&\leq& c_m\left(\mE\|X_0\|^2_\mH+\int^T_0\mE\left(\xi(s)
+\eta^{\frac{q_1}{q_1-1}}_1(s)+\eta^{\frac{q_2}{q_2-1}}_2(s)\right)\dif s\right)<+\infty,
\de
where $c_m>0$ is independent of $n$, and $m\in\mN$ is fixed in the next two steps.

{\it (Step 3)}

By the reflexivities of Banach spaces $\mK^m$, one may find a common subsequence $n_k$
(denoted by $k$ for simplicity) and
$\tilde X^m\in \mK^m_{2,1}\cap\mK^m_{2,2}$, $\bar X^m\in\mK^m_4\cap\mK^m_5$,
$Y^m_i\in\mK^m_{1,i}$, $i=1,2$, $Z^m\in\mK^m_3$
and $X^m_\infty\in L^2(\Omega,\cF_{\theta_m}, P;\mH)$ such that as $k\rightarrow\infty$
\be
X_k&\rightarrow& \tilde X^m\mbox{ weakly in $\mK^m_{2,1}$ and $\mK^m_{2,2}$},\label{e7}\\
X_k&\rightarrow& \bar X^m\mbox{ weakly in $\mK^m_4$ and $\mK^m_5$},\label{e8}\\
A_i(\cdot,X_k)=:Y_{k,i}&\rightarrow& Y^m_i\mbox{ weakly in $\mK^m_{1,i}$},\ \ i=1,2,\label{e9}\\
B(\cdot,X_k)\tilde\Pi_k=:Z_k&\rightarrow& Z^m\mbox{ weakly in $\mK^m_3$},\label{e10}\\
X_k(\theta_m)&\rightarrow& X^m_\infty\mbox{ weakly in $L^2(\Omega,\cF_{\theta_m}, P;\mH)$}.\label{e11}
\ee
Clearly, $\tilde X^m$,  $\bar X^m$, $Y^m_1$, $Y^m_2$ and $Z^m$ are $\cM$-measurable.
First of all, by (iv) of Lemma \ref{l1}, we have
\ce
\tilde X^m(t,\om)=\bar X^m(t,\om),\mbox{ for $\mu^m$-almost all $(t,\om)$}.
\de
Secondly, define
\be
X^m(t):=X_0+\int^{t\wedge\theta_m}_0(Y^m_1(s)+Y^m_2(s))\dif s+\int^{t\wedge\theta_m}_0Z^m(s)\dif W(s),
\label{PP1}
\ee
then $t\mapsto X^m(t,\om)$ is continuous in $\mH$ a.s. and (cf. \cite{Ro})
$$
\mE\left(\sup_{t\in[0,\theta_m]}\|X^m(t)\|^2_\mH\right)<+\infty.
$$
Moreover, we also have
\be
X^m(t,\om)=\bar X^m(t,\om)\quad\mbox{ for $\mu^m$-almost all $(t,\om)$}.\label{e12}
\ee
Indeed, let $\zeta(t)$ be any $\mH$-valued bounded and measurable process on $(\Om,\cF,P)$.
By (\ref{e6}) we have for any $k\geq n$
\ce
\mE\left(\int^{\theta_m}_0\<\Pi_n\zeta(t), X_k(t)\>_\mH\dif t\right)&=&
\mE\left(\int^{\theta_m}_0\<\Pi_n\zeta(t), \Pi_kX_0\>_\mH\dif t\right)\\
&&+\mE\left(\int^{\theta_m}_0\int^t_0[\Pi_n\zeta(t), Y_{k,1}(s)+Y_{k,2}(s)]_\mX\dif s\dif t\right)\\
&&+\mE\left(\int^{\theta_m}_0\<\Pi_n\zeta(t), J(Z_k)(t)\>_\mH\dif t\right),
\de
where $J$ is defined by (\ref{JJ}), and we have used that
$$
\<\Pi_n\zeta(t), J(\Pi_kZ_k)(t)\>_\mH=\<\Pi_n\zeta(t), \Pi_kJ(Z_k)(t)\>_\mH=\<\Pi_n\zeta(t), J(Z_k)(t)\>_\mH.
$$
Taking  limits for $k\rightarrow\infty$,
and by Fubini's theorem, (\ref{e8}) (\ref{e9}) (\ref{e10}) and (iii), (v) of Lemma \ref{l1} we obtain
\ce
\mE\left(\int^{\theta_m}_0\<\Pi_n\zeta(t), \bar X^m(t)\>_\mH\dif t\right)=
\mE\left(\int^{\theta_m}_0\<\Pi_n\zeta(t), X^m(t)\>_\mH\dif t\right).
\de
Letting $n\rightarrow\infty$ then shows (\ref{e12}) by the arbitrariness of $\zeta$.
Using the same method, by (\ref{e9}) (\ref{e10}) and (\ref{e11}) we also have
\be
X^m(\theta_m(\om),\om)=X^m_\infty(\om)\quad\mbox{\ \  for $P$-almost all $\om\in\Om$}.\label{e111}
\ee
In the following we shall not distinguish $X^m$, $\tilde X^m$ and $\bar X^m$.

{\it (Step 4)}

Our task in this step is to show by the standard monotone argument that for $\mu^m$-almost all
$(t,\om)\in[0,T]\times\Om$
\ce
(A_1+A_2)(t,\om, X^m(t, \om))&=&(Y^m_1+Y^m_2)(t,\om)=:Y^m(t,\om),\\
 B(t,\om, X^m(t, \om))&=&Z^m(t,\om).
\de
By Ito's formula and $\mathbf{(H2)}$(with $\lambda_0=0$), from (\ref{e6}) we have for any
$\Phi\in \mK^m_{2,1}\cap\mK^m_{2,2}\cap \mK^m_5$
\be
&&\|X_k(\theta_m)\|^2_\mH-\|\Pi_kX_0\|^2_\mH-2M(\theta_m)\no\\
&=&\int^{\theta_m}_0\left(2[X_k(s), A(s,X_k(s))]_\mX+\|\Pi_kB(s,X_k(s))\tilde\Pi_k\|^2_{L_2(\mU,\mH)}\right)
\dif s\no\\
&\leq&\int^{\theta_m}_0\left(2[X_k(s), A(s,X_k(s))]_\mX+\|B(s,X_k(s))\|^2_{L_2(\mU,\mH)}\right)\dif s\no\\
&\leq&\int^{\theta_m}_0\Big(2[X_k(s), A(s,\Phi(s))]_\mX
+2[\Phi(s), A(s,X_k(s))-A(s,\Phi(s))]_\mX\no\\
&&-\|B(s,\Phi(s))\|^2_{L_2(\mU,\mH)}+2\<B(s,X_k(s)),B(s,\Phi(s))\>_{L_2(\mU,\mH)}\Big)\dif s,\label{e13}
\ee
where $M(t)$ is a continuous martingale defined by
\ce
M(t):=\int^{t\wedge\theta_m}_0\<X_k(s), B(s,X_k(s))\tilde\Pi_k\dif W(s)\>_\mH.
\de
Since $\Phi\in \mK^m_{2,1}\cap\mK^m_{2,2}\cap \mK^m_5$, we have $B(\cdot,\Phi(\cdot))\in\mK^m_3$ by Remark \ref{r1}.
Firstly taking expectations  for (\ref{e13}), and then taking limits for $k\rightarrow\infty$,
we find by (\ref{e7}) (\ref{e9}) (\ref{e10})
and (iii) (iv) of Lemma \ref{l1} that
\ce
&&\liminf_{k\rightarrow\infty}\mE\|X_k(\theta_m)\|^2_\mH-\mE\|X_0\|^2_\mH\\
&\leq&\mE\Bigg(\int^{\theta_m}_0\Big(2[X^m(s), A(s,\Phi(s))]_\mX
+2[\Phi(s), Y^m(s)-A(s,\Phi(s))]_\mX\\
&&-\|B(s,\Phi(s))\|^2_{L_2(\mU,\mH)}+2\<Z^m(s),B(s,\Phi(s))\>_{L_2(\mU,\mH)}\Big)\dif s\Bigg).
\de
On the other hand, from (\ref{PP1}), noting that by It\^o's formula again
$$
\mE\|X^m(\theta_m)\|^2_\mH-\mE\|X_0\|^2_\mH=\mE\left(\int^{\theta_m}_0
\left(2[X^m(s),Y^m(s)]_\mX+\|Z^m(s)\|^2_{L_2(\mU,\mH)}\right)\dif s\right),
$$
and by (\ref{e11}) and (\ref{e111})
$$
\mE\|X^m(\theta_m)\|^2_\mH\leq\liminf_{k\rightarrow\infty}\mE\|X_k(\theta_m)\|^2_\mH,
$$
we finally arrive at
\ce
&&\mE\left(\int^{\theta_m}_02[X^m(s)-\Phi(s), Y^m(s)-A(s,\Phi(s))]_\mX\dif s\right)\no\\
&&+\mE\left(\int^{\theta_m}_0\|B(s,\Phi(s))-Z^m(s)\|^2_{L_2(\mU,\mH)}\dif s\right)\leq0.
\de
Letting $\Phi=X^m$ in the above inequality, we obtain that $Z^m=B(\cdot,X^m)$.
By Lemma \ref{Le2} we also have $Y^m=A(\cdot, X^m)$.

{\it (Step 5)}

For  $m\geq l$, since $\theta_m(\om)\geq\theta_l(\om)$ a.s.,
both $X^m(\cdot,\om)$ and $X^l(\cdot,\om)$ solve the following equation
\ce
\dif X(t)=A(t,X(t))1_{\{t\leq\theta_l\}}\dif t+B(t,X(t))1_{\{t\leq\theta_l\}}\dif W(s),\quad  X(0)=X_0.
\de
The uniqueness of solutions gives that for almost all $\om$
\ce
X^m(t,\om)=X^l(t,\om),\quad t\leq\theta_l(\om).
\de
Thus, noting that $\theta_m(\om)\uparrow T$ a.s. as $m\uparrow\infty$, we may define a continuous  $\cF_t$-adapted
$\mH$-valued process for all $t\in(0,T)$ by
\ce
X(t,\om):=X^{m}(t,\om)\quad\mbox { if $t<\theta_m(\om)$},
\de
Clearly, it is a solution of Eq.(\ref{Eq1}) in the sense of Definition \ref{d2}.
The proof is complete.
\end{proof}
\br
Since $q_1,q_2\geq 2$ is only used in Step 1,
if $\lambda_0=0$, Theorem \ref{th1} still holds for any $q_1,q_2>1$.
\er

\section{Backward Stochastic Evolution Equations}

In this section, we consider the following type of backward stochastic evolution equation:
\be
\left\{
\begin{array}{ll}
\dif X(t)=-A(t,X(t))\dif t-C(t,X(t),Z(t))\dif t+ Z(t)\dif W(t),\\
X(T)=X_T\in\cF_T/\cB(\mH),
\end{array}
\right.\label{Beq}
\ee
where
\ce
A_i: [0,T]\times\Om\times\mX_i\to\mX^*_i&\in& \cM\times\cB(\mX_i)/\cB(\mX^*_i),\ \ i=1,2,\\
C: [0,T]\times\Om\times\mH\times L_2(\mU,\mH)\to\mH&\in&\cM\times\cB(\mH)\times\cB(L_2(\mU,\mH))/\cB(\mH).
\de

We assume that
\begin{enumerate}[(\bf $\mathbf{HB}$1)]
\item $X_T\in L^2(\Om,\cF_T, P; \mH)$ and $(A,0)$ satisfies
$\sH(\lambda_0,\lambda_1,\lambda_2,\lambda_3,\xi,\eta_1,\eta_2, q_1,q_2)$,
where $\l_i$, $i=0,1,2,3$ are positive constants, $q_i\geq 2$,
$0<\xi\in L^1(\fA)$ and
$$
\mE\left(\int^T_0|\eta_1(s)+\eta_2(s)|^2\dif s\right)^{(q_1\vee q_2)/2}<+\infty.
$$

\item There exist a $c_1>0$ and an increasing  concave function $\rho$ satisfying
(\ref{rho}) such that for all $(t,\om)\in[0,T]\times\Om$,
$x,x'\in\mH$ and $z,z'\in L_2(\mU,\mH)$
\ce
\|C(t,\om,x,z)-C(t,\om,x',z')\|^2_\mH\leq c_1\left(\rho(\|x-x'\|^2_\mH)+\|z-z'\|^2_{L_2(\mU,\mH)}\right).
\de

\item There exist a $c_2>0$ and a $0<\zeta\in L^2(\fA)$ such that for all $(t,\om)\in[0,T]\times\Om$, $x\in\mH$ and $z\in L_2(\mU,\mH)$
$$
\|C(t,\om,x,z)\|_\mH\leq \zeta(t,\om)+c_2\left(\|x\|_\mH+\|z\|_{L_2(\mU,\mH)}\right).
$$
\end{enumerate}

Recalling Remark \ref{r3}, we give the following definition.
\bd\label{d6}
A pair of measurable $\cF_t$-adapted processes $(X,Z)$ is called a solution of Eq.(\ref{Beq}) if
\begin{enumerate}[(i)]
\item $X\in\mK_{2,1}\cap\mK_{2,2}\cap\mK_4$ and $Z\in\mK_3$, $X(0)\in L^2(\Om,\cF_0,P;\mH)$.

\item For almost all $\om$, $t\mapsto X(t,\om)$ is continuous in $\mH$ and $X(T)=X_T$ a.s..

\item $(X,Z)$ satisfies that  in $\mX^*$, for almost all $\om$ and all $t\in[0,T]$
$$
X(t)=X_T+\int^T_tA(s,X(s))\dif s+\int^T_tC(s,X(s),Z(s))\dif s-\int^T_t Z(s)\dif W(s).
$$
\end{enumerate}
\ed

As in Step 1 of Theorem \ref{th1}, let $\gamma(t):=e^{\lambda_0 t/2}$ and define
\ce
(\tilde X(t,\omega),\tilde Z(t,\omega))&:=&(\gamma(t)\cdot X(t,\omega),\gamma(t)\cdot Z(t,\omega)),\\
\tilde A_i(t,\omega,x)&:=&\gamma(t)\cdot A_i(t,\omega,\gamma^{-1}(t)\cdot x)-\lambda_0  \cdot x/2,\ \ i=1,2,\\
\tilde C(t,\omega,x,z)&:=&\gamma(t)\cdot C(t,\omega,\gamma^{-1}(t)\cdot x,\gamma^{-1}(t)\cdot z).
\de
Then we can assume $\l_0=0$ in $\mathbf{(HB1)}$ in the following.

We have the following uniqueness result.
\bt
Assume that $\mathbf{(HB1)}$, $\mathbf{(HB2)}$ and $\mathbf{(HB3)}$ hold. Let $(X,Z)$ and $(\tilde X,\tilde Z)$
be two solutions of Eq.(\ref{Beq}) with the same terminal values $X_T$. Then
for $(\dif t\times \dif P)$-almost all $(t,\om)\in[0,T]\times\Om$
$$
X(t,\om)=\tilde X(t,\om),\quad Z(t,\om)=\tilde Z(t,\om).
$$
\et
\begin{proof}
Set $Y(t):=X(t)-\tilde X(t)$. By It\^o's formula(Theorem \ref{Ito}), we have
\ce
&&\|Y(t)\|^2_\mH+\int^T_t\|Z(s)-\tilde Z(s)\|^2_{L_2(\mU,\mH)}\dif s\\
&=&2\int^T_t[Y(s), A(s,X(s))-A(s,\tilde X(s))]_\mX\dif s\\
&&+2\int^T_t\<Y(s), C(s,X(s),Z(s))-C(s,\tilde X(s),\tilde Z(s))\>_\mH\dif s\\
&&-2\int^T_t\<Y(s),(Z(s)-\tilde Z(s))\dif W(s)\>_\mH.
\de
Taking expectations, by $\mathbf{(H2)}$(with $\lambda_0=0$),
$\mathbf{(HB2)}$  and Young's inequality (\ref{Young}) we have
\ce
&&\mE\|Y(t)\|^2_\mH+\int^T_t\mE\|Z(s)-\tilde Z(s)\|^2_{L_2(\mU,\mH)}\dif s\\
&\leq&c_0\mE\left(\int^T_t\|Y(s)\|_\mH \left(\rho(\|Y(s)\|^2_\mH)
+\|Z(s)-\tilde Z(s)\|^2_{L_2(\mU,\mH)}\right)^{1/2}\dif s\right)\\
&\leq&c_0\int^T_t\mE\|Y(s)\|^2_\mH\dif s
+\frac{1}{2}\int^T_t\mE\rho(\|Y(s)\|^2_\mH)\dif s\\
&&+\frac{1}{2}\int^T_t\mE\|Z(s)-\tilde Z(s)\|^2_{L_2(\mU,\mH)} \dif s.
\de
Hence, by Jensen's inequality
\ce
&&\mE\|Y(t)\|^2_\mH+\frac{1}{2}\int^T_t\mE\|Z(s)-\tilde Z(s)\|^2_{L_2(\mU,\mH)}\dif s\\
&\leq& c_0\int^T_t\mE\|Y(s)\|^2_\mH\dif s+\frac{1}{2}\int^T_t\rho\left(\mE\|Y(s)\|^2_\mH\right)\dif s.
\de
The uniqueness follows from Lemma \ref{le1}.
\end{proof}

The following finite dimensional result was proved in \cite{Br-Ca}.
For completeness, we give a different proof by Yosida's approximation.
\bl\label{l0}
Assume that $\mX=\mH=\mX^*=\mU=\mR^d$ and $C=0$,
 and $(A,0)$ satisfies
$\sH(0,\lambda_1,0,\lambda_3,\xi,\eta,0, q,0)$,
where $\l_1,\l_3$ are positive constants, $q\geq 2$,
$0<\xi\in L^1(\fA)$ and
$$
\mE\left(\int^T_0|\eta(s)|^2\dif s\right)^{q/2}<+\infty.
$$
Then  for any $X_T\in L^q(\Omega,\cF_T,P;\mR^d)$,
 there exists a unique solution to Eq.(\ref{Beq}) in the sense of
Definition \ref{d6}. Moreover,
\be
&&\mE\left(\sup_{t\in[0,T]}\|X(t)\|^{q}_{\mR^d}\right)
+\mE\left(\int^T_0\|Z(s)\|^2_{L_2(\mR^d,\mR^d)}\dif s\right)^{q/2}\no\\
&\leq& c_0\cdot\mE\|X_T\|^q_{\mR^d}
+c_0\cdot\mE\left(\int^T_0|\eta(s)|^2\dif s\right)^{q/2},\label{es03}
\ee
where $c_0$ only depends on $q,T$ and $\lambda_1$.
\el
\begin{proof}
For every $(t,\om)\in[0,T]\times\Om$, note that
$x\mapsto A(t,\om, x)$ is a continuous monotone function on $\mR^d$.
Let $A_\e(t,\om,\cdot), \e>0$ be the Yosida approximation of $A(t,\om, \cdot)$, i.e.:
\ce
A_\e(t,\om,x)&:=&\e^{-1}(J_\e(t,\om,x)-x)=A(t,\om,J_\e(t,\om,x)),\\
J_\e(t,\om,x)&:=&(I-\e A(t,\om,\cdot))^{-1}(x),
\de
then $x\mapsto J_\e(t,\om,x)$ is a homeomorphism on $\mR^d$ for each $(t,\om)$ and
for any $x,y\in\mR^d$(cf. \cite{Ba} \cite{Br})
\begin{enumerate}[(I)]
\item $\<x-y, A_\e(t,\om,x)-A_\e(t,\om,y)\>_{\mR^d}\leq 0$,

\item $\|A_\e(t,\om,x)-A_\e(t,\om,y)\|_{\mR^d}\leq \e^{-1}\|x-y\|_{\mR^d}$,

\item $\|A_\e(t,\om,x)\|_{\mR^d}\leq \|A(t,\om,x)\|_{\mR^d}$,

\item $\lim_{\e\downarrow 0}\|A_\e(t,\om,x)-A(t,\om,x)\|_{\mR^d}=0$.
\end{enumerate}

By Lemma \ref{l4}, $J_\e$ and $A_\e$ are progressively measurable.
From (I), (III) and $\mathbf{(H4)}$, we have for any $x\in\mR^d$
\be
\<x, A_\e(t,\om,x)\>_{\mR^d}\leq \|x\|_{\mR^d}\cdot\|A(t,\om,0)\|_{\mR^d}
\leq \eta(t,\om)\cdot\l_1^{\frac{1}{q}}\cdot\|x\|_{\mR^d}.\label{Op2}
\ee
Let $(X_\e,Z_\e)$ be the unique $\cF_t$-adapted
solution of the following backward stochastic differential equation(cf. \cite{Pa-Pe})
\be
X_\e(t)=X_T+\int^T_tA_\e(s,X_\e(s))\dif s-\int^T_tZ_\e(s)\dif W(s).\label{es8}
\ee
By It\^o's formula, we have
\be
&&\|X_\e(t)\|^2_{\mR^d}+\int^T_t\|Z_\e(s)\|^2_{L_2(\mR^d,\mR^d)}\dif s\no\\
&=&\|X_T\|^{2}_{\mR^d}+2\int^T_t\<X_\e(s),A_\e(s,X_\e(s))\>_{\mR^d}\dif s\no\\
&&-2\int^T_t\<X_\e(s),Z_\e(s)\dif W(s)\>_{\mR^d}\label{es7}
\ee
and further
\be
&&e^{ t}\|X_\e(t)\|^2_{\mR^d}+\int^T_te^{ s}\left(\|X_\e(s)\|^2_{\mR^d}
+\|Z_\e(s)\|^2_{L_2(\mR^d,\mR^d)}\right)\dif s\no\\
&=&e^{ T}\|X_T\|^{2}_{\mR^d}+2\int^T_te^{ s}\<X_\e(s),A_\e(s,X_\e(s))\>_{\mR^d}\dif s\no\\
&&-2\int^T_te^{ s}\<X_\e(s),Z_\e(s)\dif W(s)\>_{\mR^d}\no\\
&\leq&e^{ T}\|X_T\|^{2}_{\mR^d}+
\int^T_te^{ s}\|X_\e(s)\|^2_{\mR^d}\dif s+c_0\int^T_te^{ s}|\eta(s)|^2\dif s\no\\
&&-2\int^T_te^{ s}\<X_\e(s),Z_\e(s)\dif W(s)\>_{\mR^d},\label{es00}
\ee
where the second step is due to (\ref{Op2}) and Young's inequality (\ref{Young}).

Taking conditional expectations for both sides of (\ref{es00}) with respect to $\cF_t$,
we find
\ce
e^t\|X_\e(t)\|^2_{\mR^d}&\leq& e^T\mE^{\cF_t}\|X_T\|^2_{\mR^d}+
c_0\cdot\mE^{\cF_t}\int^T_te^s |\eta(s)|^2\dif s\\
&\leq& e^T\mE^{\cF_t}\|X_T\|^2_{\mR^d}+
c_0\cdot e^T\cdot\mE^{\cF_t}\int^T_0|\eta(s)|^2\dif s.
\de
Hence, by Doob's maximal inequality(cf. \cite{Re-Yo}),  we have for $q>2$
\be
\mE\left(\sup_{t\in[0,T]}\|X_\e(t)\|^{q}_{\mR^d}\right)\leq c_0\cdot\mE\|X_T\|^q_{\mR^d}
+c_0\cdot\mE\left(\int^T_0|\eta(s)|^2\dif s\right)^{q/2}.\label{es02}
\ee
Hereafter $c_0$ only depends on $q,T$ and $\lambda_1$.

Noting that by BDG's inequality and Young's inequality (\ref{Young})
\ce
&&\mE\left|\int^T_0e^s\<X_\e(s),Z_\e(s)\dif W(s)\>_{\mR^d}\right|^{q/2}\\&\leq&
c_0\mE\left(\int^T_0e^{2s}\|X_\e(s)\|^2_{\mR^d}\cdot\|Z_\e(s)\|^2_{L_2(\mR^d,\mR^d)}\dif s\right)^{q/4}\\
&\leq&c_0\mE\left(\sup_{t\in[0,T]}\|X_\e(t)\|^{q}_{\mR^d}\right)
+\frac{1}{2}\mE\left(\int^T_0\|Z_\e(s)\|^2_{L_2(\mR^d,\mR^d)}\dif s\right)^{q/2},
\de
we also have from (\ref{es00})
\be
\mE\left(\int^T_0\|Z_\e(s)\|^2_{L_2(\mR^d,\mR^d)}\dif s\right)^{q/2}\leq c_0\cdot\mE\|X_T\|^q_{\mR^d}
+c_0\cdot\mE\left(\int^T_0|\eta(s)|^2\dif s\right)^{q/2}.\label{es01}
\ee

For $q=2$, from (\ref{es00}) and the above proof, it is easy to see that
\ce
&&\mE\left(\sup_{t\in[0,T]}\|X_\e(t)\|^2_{\mR^d}\right)
+\int^T_0\mE\|Z_\e(s)\|^2_{L_2(\mR^d,\mR^d)}\dif s\\
&\leq& c_0\cdot\mE\|X_T\|^2_{\mR^d}
+c_0\cdot\mE\int^T_0|\eta(s)|^2\dif s.
\de

Moreover, by (III), $\mathbf{(H4)}$  and (\ref{es02})
\ce
\int^T_0\|A_\e(s,X_\e(s))\|^{\frac{q}{q-1}}_{\mR^d}\dif s
&\leq&\int^T_0\|A(s,X_\e(s))\|^{\frac{q}{q-1}}_{\mR^d}\dif s\\
&\leq&c_0\int^T_0\|X_\e(s)\|^q_{\mR^d}\dif s+c_0\int^T_0\mE\left(\eta^{\frac{q}{q-1}}(s)\right)\dif s
\leq c'_0.
\de
Therefore, there exists a subsequence $\e_n\downarrow 0$ and $(X,Y,Z,X_0)$ such that
\ce
X_{\e_n}&\rightarrow& X\mbox{ weakly in $\mK_{2,1}$},\\
A_{\e_n}(\cdot,X_{\e_n}(\cdot))&\rightarrow& Y
\mbox{ weakly in $\mK_{1,1}$},\\
Z_{\e_n}&\rightarrow& Z\mbox{ weakly in $\mK_3$},\\
X_{\e_n}(0)&\rightarrow& X_0\mbox{ weakly in $L^2(\Om,\cF_0, P;\mR^d)$}
\de
as $n\rightarrow\infty$. By (\ref{es02}) and (\ref{es01}), we get (\ref{es03}).

Set
$$
\tilde X(t):=X_T+\int^T_tY(s)\dif s-\int^T_tZ(s)\dif W(s).
$$
By taking  weak limits for (\ref{es8}), we deduce that $\tilde X(0)=X_0$ a.s. and
\ce
X(t,\om)=\tilde X(t,\om)\mbox{ for almost all $(t,\om)\in[0,T]\times\Om$}.
\de

It remains to show that $Y(s)=A(s,X(s))$. For any $\Phi\in \mK_{2,1}$,
by (III) (IV) and the dominated convergence theorem we have
\be
&&\lim_{n\rightarrow\infty}\mE\left(\int^T_0\<X_{\e_n}(s)-\Phi(s),
A_{\e_n}(s,\Phi(s))-A(s,\Phi(s))\>_{\mR^d}\dif s\right)\no\\
&\leq&\lim_{n\rightarrow\infty}
\left(\|A_{\e_n}(\cdot,\Phi(\cdot))-A(\cdot,\Phi(\cdot))\|^{\frac{q}{q-1}}_{\mK_{1,1}}
\cdot\|X_{\e_n}-\Phi\|_{L^q(\fA)}\right)
=0.\label{es10}
\ee
On the other hand, we have by (\ref{es7})
\be
&&2\liminf_{n\rightarrow\infty}\mE\left(\int^T_0
\<X_{\e_n}(s),A_{\e_n}(s,X_{\e_n}(s))\>_{\mR^d}\dif s\right)\label{es11}\\
&\geq&\mE\|X_0\|^{2}_{\mR^d}-\mE\|X_T\|^{2}_{\mR^d}+\int^T_0\mE\|Z(s)\|^2_{L_2(\mR^d,\mR^d)}\dif s\no\\
&=&2\mE\left(\int^T_0\<X(s),Y(s))\>_{\mR^d}\dif s\right).\no
\ee
Combining  (\ref{es10}) and (\ref{es11}), we have by (I)
\ce
&&\mE\left(\int^T_0\<X(s)-\Phi(s),Y(s)-A(s,\Phi(s)))\>_{\mR^d}\dif s\right)\\
&\leq&\liminf_{n\rightarrow\infty}
\mE\left(\int^T_0\<X_{\e_n}(s)-\Phi(s),A_{\e_n}(s,X_{\e_n}(s))-A_{\e_n}(s,\Phi(s))\>_{\mR^d}\dif s\right)\leq 0,
\de
which implies that $Y=A(\cdot,X)$ by Lemma \ref{Le2}. The proof is complete.
\end{proof}
\br
When $q>2$, it suffices to require that
\ce
\mE\left(\int^T_0|\eta(s)|\dif s\right)^{q}<+\infty.
\de
In fact, taking conditional expectations for both sides of (\ref{es7})
with respect to $\cF_t$, and by (\ref{Op2}) and Young's inequality (\ref{Young})
we find for any $\delta>0$
\ce
\|X_\e(t)\|^2_{\mR^d}&\leq& \mE^{\cF_t}\|X_T\|^2_{\mR^d}+q\mE^{\cF_t}\left(\int^T_t\|X_\e(s)\|_{\mR^d}
\cdot\eta(s)\cdot\l_1^{\frac{1}{q}}\dif s\right)\\
&\leq& \mE^{\cF_t}\|X_T\|^2_{\mR^d}+\delta\cdot\mE^{\cF_t}\left(\sup_{s\in[0,T]}\|X_\e(s)\|^2_{\mR^d}\right)\\
&&+c_\delta\cdot\mE^{\cF_t}\left(\int^T_0|\eta(s)|\dif s\right)^2.
\de
Hence, by Doob's maximal inequality  we have for $q>2$
\ce
\mE\left(\sup_{t\in[0,T]}\|X_\e(t)\|^{q}_{\mR^d}\right)&\leq& c_0\cdot\mE\|X_T\|^q_{\mR^d}
+\delta\cdot c_q\cdot\mE\left(\sup_{s\in[0,T]}\|X_\e(s)\|^q_{\mR^d}\right)\\
&&+c_\delta\cdot\mE\left(\int^T_0|\eta(s)|\dif s\right)^{q}.
\de
Letting $\delta$ be sufficiently small, we get
\ce
\mE\left(\sup_{t\in[0,T]}\|X_\e(t)\|^{q}_{\mR^d}\right)\leq c_0\cdot\mE\|X_T\|^q_{\mR^d}
+c_0\cdot\mE\left(\int^T_0|\eta(s)|\dif s\right)^{q},
\de
where $c_0$ only depends on $q,\lambda_1$ and $T$.
\er
We now prove the following infinite dimensional version.
\bl\label{l2}
Assume that $C(\cdot,x,z)=C(\cdot)\in L^2(\fA_a;\mH)$ is independent of $x$ and $z$, and $\mathbf{(HB1)}$ holds.
Then there exists a unique solution to Eq.(\ref{Beq}) in the sense of
Definition \ref{d6}.
\el
\begin{proof}
We use Galerkin's approximation to prove the existence as in the proof of Theorem \ref{th1}.
For $n\in\mN$, let $(X_n,Z_n)$ solve the following finite dimensional backward stochastic differential equation
(Lemma \ref{l0})
\ce
X_n(t)=X^n_T+\int^T_t\Pi_nA(s,X_n(s))\dif s+\int^T_tC_n(s)\dif s-\int^T_tZ_n(s)\tilde\Pi_n\dif W(s),
\de
where $\Pi_n$ and $\tilde\Pi_n$ are same as in Theorem \ref{th1}, and
\ce
X^n_T&:=&\Pi_n X_T\cdot 1_{\{\|\Pi_n X_T\|_\mH\leq n\}}\\
C_n(s)&:=&\Pi_n C(s)\cdot 1_{\{\|\Pi_n C(s)\|_\mH\leq n\}}.
\de
It is easy to see that for each $n$ and $s$
$$
\|X^n_T\|_\mH\leq \|X_T\|_\mH,\quad \|C_n(s)\|_\mH\leq\|C(s)\|_\mH
$$
and
\be
\lim_{n\rightarrow\infty}\mE\|X^n_T-X_T\|^2_\mH&=&0,\label{Op4}\\
\lim_{n\rightarrow\infty}\int^T_0\mE\|C_n(s)-C(s)\|^2_\mH\dif s&=&0.\label{Op5}
\ee

By It\^o's formula and $\mathbf{(H3)}$, we have
\be
&&\mE\|X_n(t)\|^2_\mH+\int^T_t\mE\|Z_n(s)\|^2_{L_2(\mU,\mH)}\dif s\label{Op3}\\
&=&\mE\|X^n_T\|^2_\mH+\int^T_t\mE[X_n(s),A(s,X_n(s))]_\mX\dif s
+2\int^T_t\mE\<X_n(s),C_n(s)\>_\mH\dif s\no\\
&\leq&\mE\|X_T\|^2_\mH+\int^T_0\mE(2\xi(s)+\|C(s)\|^2_\mH)\dif s\no\\
&&+\int^T_t\mE\left(-\sum_{i=1,2}\l_i\cdot\|X_n(s)\|^{q_i}_{\mX_i}
+(\l_3+1)\cdot\|X_n(s)\|^2_\mH\right)\dif s.\label{Es1}
\ee
By Gronwall's inequality we have
\ce
\mE\|X_n(t)\|^2_\mH\leq c_0\left(\mE\|X_T\|^2_\mH+\int^T_0\mE(\xi(s)+\|C(s)\|^2_\mH)\dif s\right).
\de
Hence, from (\ref{Es1}) and $\mathbf{(H4)}$ we get
\ce
\mE\|X_n(0)\|^2_\mH+\|X_n\|^2_{\mK_4}+\|Z_n\|^2_{\mK_3}
+\sum_{i=1,2}\Big(\|X_n\|^{q_i}_{\mK_{2,i}}+\|A_i(\cdot,X_n)\|^{\frac{q_i}{q_i-1}}_{\mK_{1,i}}\Big)\leq c_0.
\de
Hereafter, the constant $c_0$ is independent of $n$.

By the reflexivities  of Banach spaces $\mK$, one may find a subsequence
$n_k$ (denoted by $k$ for simplicity) and $\tilde X\in \mK_{2,1}\cap\mK_{2,2}\cap  \mK_4$, $Y_i\in\mK_{1,i}$,
$i=1,2$ and $Z\in\mK_3$ such that
\ce
X_k&\to& \tilde X \mbox{ weakly in $\mK_{2,1},\mK_{2,2}$ and $\mK_4$},\\
A_i(\cdot,X_k)=:Y_{k,i}&\to& Y_i\mbox{ weakly in $\mK_{1,i}$},\ \ i=1,2,\\
Z_k&\to& Z\mbox{ weakly in $\mK_3$},\\
X_k(0)&\to& X_0\mbox{ weakly in $L^2(\Om,\cF_0,P; \mH)$}.
\de
Define $Y=Y_1+Y_2\in\mY\subset\mX^*$ and
$$
X(t):=X_T+\int^T_tY(s)\dif s+\int^T_tC(s)\dif s-\int^T_tZ(s)\dif W(s).
$$
Then, similar to  Step 3 of Theorem \ref{th1} one may prove that
\ce
\tilde X(t,\om)&=&X(t,\om)\quad\mbox{ for ($\dif t\times \dif P$)-almost all $(t,\om)$},\\
X(0)&=&X_0 \quad a.s..
\de

We now show that
\be
A(t, X(t,\om))=Y(t,\om)\quad\mbox{ for ($\dif t\times \dif P$)-almost all $(t,\om)$}.\label{e15}
\ee
By (\ref{Op3}) and $\mathbf{(H2)}$(with $\lambda_0=0$), we have for any $\Phi\in\mK_2$
\ce
&&\mE\|X_k(0)\|^2_\mH+\int^T_0\mE\|Z_k(s)\|^2_{L_2(\mU,\mH)}\dif s\\
&\leq&\mE\|X^k_T\|^2_\mH+2\int^T_0\mE[\Phi(s),A(s,X_k(s))-A(s,\Phi(s))]_\mX\dif s\\
&&+2\int^T_0\mE[X_k(s),A(s,\Phi(s))]_\mX\dif s+2\int^T_0\mE\<X_k(s),C_k(s)\>_\mH\dif s.
\de
Taking  limits for $k\rightarrow\infty$, we find by (\ref{Op4}) (\ref{Op5})
\ce
&&\liminf_{k\rightarrow\infty}\mE\|X_k(0)\|^2_\mH
+\liminf_{k\rightarrow\infty}\int^T_0\mE\|Z_k(s)\|^2_{L_2(\mU,\mH)}\dif s\\
&\leq&\mE\|X_T\|^2_\mH+2\int^T_0\mE[\Phi(s),Y(s)-A(s,\Phi(s))]_\mX\dif s\\
&&+2\int^T_0\mE[X(s),A(s,\Phi(s))]_\mX\dif s+2\int^T_0\mE\<X(s),C(s)\>_\mH\dif s.
\de
On the other hand, noting that
\ce
&&\mE\|X_0\|^2_\mH+\int^T_0\mE\|Z(s)\|^2_{L_2(\mU,\mH)}\dif s\\
&=&\mE\|X_T\|^2_\mH+2\int^T_0\mE[X(s),Y(s)]_\mX\dif s+2\int^T_0\mE\<X(s),C(s)\>_\mH\dif s
\de
and
\ce
\mE\|X_0\|^2_\mH&\leq&\liminf_{k\rightarrow\infty}\mE\|X_k(0)\|^2_\mH,\\
\int^T_0\mE\|Z(s)\|^2_{L_2(\mU,\mH)}\dif s&\leq&
\liminf_{k\rightarrow\infty}\int^T_0\mE\|Z_k(s)\|^2_{L_2(\mU,\mH)}\dif s,
\de
we obtain
$$
\int^T_0\mE[X(s)-\Phi(s),Y(s)-A(s,\Phi(s))]_\mX\dif s\leq 0.
$$
Hence $Y=A(\cdot,X)$ by Lemma \ref{Le2}.
The proof is complete.
\end{proof}
\bl\label{l3}
Assume that $C(t,x,z)=C(t,z)$ is independent of
$x$, and $\mathbf{(HB1)}$, $\mathbf{(HB2)}$ and $\mathbf{(HB3)}$ hold. Then there exists a unique solution
to Eq.(\ref{Beq}) in the sense of
Definition \ref{d6}.
\el
\begin{proof}
Let $Z_0(t)\equiv 0$. We consider the following Picard iteration: for $n\in\mN$,
let $(X_n,Z_n)$ solve the following equation(Lemma \ref{l2}):
$$
X_{n}(t)=X_T+\int^T_tA(s,X_{n}(s))\dif s+\int^T_tC(s,Z_{n-1}(s))\dif s-\int^T_tZ_{n}(s)\dif W(s).
$$

Set $Y_n(t):=X_{n+1}(t)-X_n(t)$.
By It\^o's formula,
$\mathbf{(H2)}$(with $\lambda_0=0$), $\mathbf{(HB2)}$ and Young's inequality, we have
\be
&&\mE\|Y_n(t)\|^2_\mH+\int^T_t\mE\|Z_{n+1}(s)-Z_{n}(s)\|^2_{L_2(\mU,\mH)}\dif s\no\\
&=&\int^T_t\mE[Y_n(s),A(s,X_{n+1}(s))-A(s,X_{n}(s))]_\mX\dif s\no\\
&&+\int^T_t\mE\<Y_n(s),C(s,Z_n(s))-C(s,Z_{n-1}(s))\>_\mH\dif s\no\\
&\leq&c_0\int^T_t\mE\|Y_n(s)\|^2_\mH\dif s+\frac{1}{2}\int^T_t\|Z_n(s)-Z_{n-1}(s)\|^2_{L_2(\mU,\mH)}\dif s.\label{es1}
\ee
Hence, for $\a:=c_0$
\ce
&&-\frac{\dif }{\dif t}\left(e^{\a t}\int^T_t\mE\|Y_n(s)\|^2_\mH\dif s\right)
+e^{\a t}\int^T_t\mE\|Z_{n+1}(s)-Z_{n}(s)\|^2_{L_2(\mU,\mH)}\dif s\\
&\leq&\frac{e^{\a t}}{2}\int^T_t\|Z_n(s)-Z_{n-1}(s)\|^2_{L_2(\mU,\mH)}\dif s=:\frac{g_n(t)}{2}.
\de
Integrating both sides from $0$ to $T$ yields that
\ce
\int^T_0\mE\|Y_n(s)\|^2_\mH\dif s+\int^T_0g_{n+1}(t)\dif t
\leq\frac{1}{2}\int^T_0g_{n}(t)\dif t\leq \frac{1}{2^n}\int^T_0g_1(t)\dif t=:\frac{c_0}{2^n}.
\de
It then follows from (\ref{es1}) that
\ce
\int^T_0\mE\|Z_{n+1}(s)-Z_{n}(s)\|^2_{L_2(\mU,\mH)}\dif s\leq\frac{c_0}{2^n}
+\frac{1}{2}\int^T_0\mE\|Z_{n}(s)-Z_{n-1}(s)\|^2_{L_2(\mU,\mH)}\dif s.
\de
Iterating this inequality gives
$$
\int^T_0\mE\|Z_{n+1}(s)-Z_{n}(s)\|^2_{L_2(\mU,\mH)}\dif s\leq\frac{nc_0}{2^n}.
$$
Therefore, there exist an $X\in\mK_4$ and a $Z\in\mK_3$ such that
\ce
\lim_{n\rightarrow\infty}\|X_n-X\|_{\mK_4}=0\mbox{ and } \lim_{n\rightarrow\infty}\|Z_n-Z\|_{\mK_3}=0.
\de
From (\ref{es1}) and the above estimates, we also have
\be
\sup_{n\in\mN}\sup_{t\in[0,T]}\mE\|X_n(t)\|^2_\mH<+\infty.\label{es3}
\ee

We now show that there exists a version $(\tilde X,\tilde Z)$ of $(X, Z)$ such that $(\tilde X,\tilde Z)$
is a solution to Eq.(\ref{Beq}) in the sense of Definition \ref{d6}.
In fact, let $(\tilde X,\tilde Z)$ solve the following equation(Lemma \ref{l2}):
$$
\tilde X(t)=X_T+\int^T_tA(s,\tilde X(s))\dif s+\int^T_tC(s,Z(s))\dif s-\int^T_t\tilde Z(s)\dif W(s).
$$
It is similar to estimate (\ref{es1}) that
\ce
&&\mE\|X_{n}(t)-\tilde X(t)\|^2_\mH+\int^T_t\mE\|Z_{n}(s)-\tilde Z(s)\|^2_{L_2(\mU,\mH)}\dif s\\
&\leq&c_0\int^T_t\mE\|X_{n}(s)-\tilde X(s)\|^2_\mH\dif s+\frac{1}{2}\int^T_t\|Z_{n-1}(s)-Z(s)\|^2_{L_2(\mU,\mH)}\dif s.
\de
Letting $g(t):=\limsup_{n\rightarrow\infty}\mE\|X_{n}(t)-\tilde X(t)\|^2_\mH$,
by (\ref{es3}) and Fatou's lemma, we have
$$
g(t)\leq c_0\int^T_tg(s)\dif s.
$$
which yields that   $g(t)=0$ by Gronwall's inequality. The proof is complete.
\end{proof}

We now prove our main result in this section.
\bt
Assume that $\mathbf{(HB1)}$, $\mathbf{(HB2)}$ and $\mathbf{(HB3)}$ hold.
Then there exists a unique solution to Eq.(\ref{Beq}) in the sense of
Definition \ref{d6}.
\et
\begin{proof}
Let $X_0(t)\equiv0$. We consider the following Picard iteration: for $n\in\mN$,
let $(X_n,Z_n)$ solve the following equation(Lemma \ref{l3})
\ce
X_{n}(t)=X_T+\int^T_tA(s,X_{n}(s))\dif s+\int^T_tC(s,X_{n-1}(s),Z_{n}(s))\dif s-\int^T_tZ_{n}(s)\dif W(s).
\de
First of all, by It\^o's formula, $\mathbf{(H2)}$(with $\lambda_0=0$),
$\mathbf{(HB3)}$ and Young's inequality, we have
\ce
&&\mE\|X_{n}(t)\|^2_\mH+\int^T_t\mE\|Z_{n}(s)\|^2_{L_2(\mU,\mH)}\dif s\\
&=&\mE\|X_T\|^2_\mH+2\int^T_t\mE[X_{n}(s),A(s,X_n(s))]_\mX\dif s\\
&&+2\int^T_t\mE\<X_{n}(s),C(s,X_{n-1}(s),Z_n(s)))\>_\mH\dif s\\
&\leq&\mE\|X_T\|^2_\mH+2\int^T_t\mE\left(\l_3\|X_{n}(s)\|^2_\mH+\xi(s)\right)\dif s\\
&&+2\int^T_t\mE\left(\|X_{n}(s)\|_\mH\Big(\zeta(s)+c_2(\|X_{n-1}(s)\|_\mH+\|Z_n(s)\|_{L_2(\mU,\mH)})\Big)\right)\dif s\\
&\leq&\mE\|X_T\|^2_\mH+c_0\int^T_0\mE\Big(\xi(s)+\zeta^2(s)\Big)\dif s+c_0\int^T_t\mE\|X_{n}(s)\|^2_\mH\dif s\\
&&+\frac{1}{4}\int^T_t\mE\Big(\|X_{n-1}(s)\|^2_\mH+\|Z_n(s)\|^2_{L_2(\mU,\mH)}\Big)\dif s.
\de
So
\be
&&\mE\|X_{n}(t)\|^2_\mH+\int^T_t\mE\|Z_n(s)\|^2_{L_2(\mU,\mH)}\dif s\no\\
&\leq& c_0+c_0\int^T_t\left(\mE\|X_{n}(s)\|^2_\mH+\mE\|X_{n-1}(s)\|^2_\mH\right)\dif s,
\label{es4}
\ee
where $c_0$ is independent of $n$.

Set
$$
g_n(t):=\max_{1\leq k\leq n}\mE\|X_k(t)\|^2_\mH.
$$
Then
$$
g_n(t)\leq c_0+c_0\int^T_tg_n(s)\dif s,
$$
which gives that by Gronwall's inequality
\be
\max_{k\in\mN}\sup_{t\in[0,T]}\mE\|X_k(t)\|^2_\mH\leq \max_{n\in\mN}\sup_{t\in[0,T]}g_n(t)<+\infty.\label{es5}
\ee

Set $Y_{n,m}(t):=X_{n}(t)-X_m(t)$ and $G_{n,m}(t):=Z_{n}(s)-Z_{m}(s)$. By It\^o's formula,
$\mathbf{(H2)}$(with $\lambda_0=0$)  and $\mathbf{(HB2)}$, we have
\ce
&&\mE\|Y_{n,m}(t)\|^2_\mH+\int^T_t\mE\|G_{n,m}(s)\|^2_{L_2(\mU,\mH)}\dif s\\
&=&2\int^T_t\mE[Y_{n,m}(s),A(s,X_n(s))-A(s,X_m(s))]_\mX\dif s\\
&&+2\int^T_t\mE\<Y_{n,m}(s),C(s,X_{n-1}(s),Z_n(s))-C(s,X_{m-1}(s),Z_m(s))\>_\mH\dif s\\
&\leq&c_0\int^T_t\mE\left(\|Y_{n,m}(s)\|_\mH\Big(\rho(\|Y_{n-1,m-1}(s)\|^2_\mH)
+\|G_{n,m}(s)\|^2_{L_2(\mU,\mH)}\Big)^{\frac{1}{2}}\right)\dif s.
\de
Using the same method as in estimating (\ref{es4}), we have
\ce
&&\mE\|Y_{n,m}(t)\|^2_\mH+\int^T_t\mE\|G_{n,m}(s)\|^2_{L_2(\mU,\mH)}\dif s\\
&\leq& c_0\int^T_t\mE\|Y_{n,m}(s)\|^2_\mH\dif s+c_0\int^T_t\mE\rho(\|Y_{n-1,m-1}(s)\|^2_\mH)\dif s.
\de
Set
\ce
g(t):=\limsup_{n,m\rightarrow\infty}\mE\|Y_{n,m}(t)\|^2_\mH.
\de
By (\ref{es5}), Fatou's lemma and Jensen's inequality, we have
\ce
g(t)\leq c_0\int^T_t(g(s)+\rho(g(s)))\dif s.
\de
So, by Lemma \ref{le1}
\ce
g(t)=0,\quad t\in[0,T].
\de
Hence
\ce
\limsup_{n,m\rightarrow\infty}\int^T_0\left(\mE\|Y_{n,m}(s)\|^2_\mH+\mE\|G_{n,m}(s)\|^2_{L_2(\mU,\mH)}\right)\dif s=0,
\de
and there exist an $X\in\mK_4$ and a $Z\in\mK_3$ such that
\ce
\lim_{n\rightarrow\infty}\|X_{n}-X\|_{\mK_4}=0\mbox{ and }\lim_{n\rightarrow\infty}\|Z_{n}-Z\|_{\mK_3}=0.
\de
Using the same method as in the proof of Lemma \ref{l3},
we can show that $(X,Z)$ solves Eq.(\ref{Beq}). The proof is thus complete.
\end{proof}

\br
In finite dimensional case,  under rather weak assumptions on $C$,
the authors \cite{Br-Ca} proved the existence and uniqueness of Eq.(\ref{Beq}).
It is interesting that the growth of $C$ in $x$ therein can be arbitrary
(not necessary polynomial growth). We remark that in our equation,
the operator $A$ may contain a polynomial growth part in $x$. However,
it seems to be difficult to extend $A$ or $C$ to be arbitrary growth in $x$ when we use the
cutoff technique as in \cite{Br-Ca}, because $A$ is a non-linear operator and
we need to take weak limits in $L^p$-space. On the other hand, if $C$ is polynomial growth in $\mH$
with respect to $x$, it will exclude the interesting case that $C$ is a Nemytskii operator.
For example, let $\varphi(r)=-|r| r$, it is not true that
$L^2(0,1)\ni x\mapsto\varphi(x)\in L^2(0,1)$, but,  $L^4(0,1)\ni x\mapsto\varphi(x)\in L^2(0,1)$.
\er

\section{Stochastic Functional Integral Evolution Equations}

Fix $S>0$.  For any $T\geq 0$, let $\mF^T_S(\mH)$ denote the space of all
 continuous functions from $[-S,T]$ to $\mH$, which
is a separable Banach space under the supremum norm
$$
\|f\|_{\mF^T_S}:=\sup_{s\in[-S,T]}\|f(s)\|_\mH.
$$
For $s\in[-S,0]$, define $\cF_s:=\cF_0$. Suppose that $X:[-S,T]\times\Omega\to\mH$ is
a continuous $\cF_t$-adapted process, we can associate it with another continuous $\mF^0_S(\mH)$-valued
and $\cF_t$-adapted process as follows:
$$
[0,T]\times\Om\ni (t,\om)\mapsto X_t(\cdot,\om):=X(t+\cdot,\om)\in \mF^0_S(\mH).
$$
In the following, we shall use the following notations:
$$
\|X_\cdot(0)\|_{\mF^t_0}:=\sup_{s\in[0,t]}\|X_s(0)\|_\mH
=\sup_{s\in[0,t]}\|X(s)\|_\mH=:\|X\|_{\mF^t_0}.
$$

Consider the following
stochastic functional integral evolution equation:
\be
X_t(0)&=&X_{0}(0)+\int^t_{0}A(s,X_s(0))\dif s+\int^t_{0}C_1(s,X_s)\dif s+\int^t_{0}\int^s_0C_2(s,r,X_r)\dif r\dif s\no\\
&&+\int^t_{0} D_1(s,X_s)\dif W_s+\int^t_{0} \int^s_0D_2(s,r,X_r)\dif W_r\dif s,\label{Eq3}
\ee
where $X_{0}$ is an $\cF_{0}$-measurable $\mF^0_S(\mH)$-valued random variable and  $A=A_1+A_2$,
\ce
A_i: [0,T]\times\Om\times\mX_i&\to&\mX^*_i,\ \ i=1,2,\\
C_1:[0,T]\times\Omega\times \mF^0_S(\mH)&\to&\mH,\\
C_2:[0,T]\times[0,T]\times\Omega\times \mF^0_S(\mH)&\to&\mH,\\
D_1:[0,T]\times\Omega\times \mF^0_S(\mH)&\to& L_2(\mU,\mH),\\
D_2:[0,T]\times[0,T]\times\Omega\times \mF^0_S(\mH)&\to& L_2(\mU,\mH)
\de
are progressively measurable, for example, for every $0\leq t\leq T$, the mapping
$(s,\om,x)\mapsto D_2(t,s,\om,x)$ is $\cM\times\cB(\mF^0_S(\mH))/\cB(L_2(\mU,\mH))$-measurable.

We make the following assumptions:
\begin{enumerate}[(\bf $\mathbf{HF}$1)]
\item $X_0\in L^2(\Omega,\cF_0,P;\mF^0_S(\mH))$ and $(A,0)$ satisfies
$\sH(\lambda_0,\lambda_1,\lambda_2,\lambda_3,\xi,\eta_1,\eta_2, q_1,q_2)$,
where $\l_i\in L^1([0,T])$, $i=0,1,2,3$ are non-random strictly  positive functions, $q_i\geq 2$,
and $0\leq\eta_i\in L^{\frac{q_i}{q_i-1}}(\fA)$, $i=1,2$, and $0\leq\xi\in L^1(\fA)$.

\item There exist a positive real function $\lambda_3\in L^1([0,T])$
and an increasing concave function $\rho$ satisfying (\ref{rho})
such that for all $(s,\om)\in[0,T]\times\Om$ and  $x,y\in{\mF^0_S}$
\ce
\|C_1(s,\om,x)-C_1(s,\om,y)\|^2_\mH&\leq&  \lambda_3(s)\cdot\rho(\|x-y\|^2_{\mF^0_S}),\\
\|D_1(s,\om,x)-D_1(s,\om,y)\|^2_{L_2(\mU,\mH)}
&\leq&  \lambda_3(s)\cdot\rho(\|x-y\|^2_{\mF^0_S}).
\de

\item There exists a positive real function $\lambda_5$ satisfying
$t\mapsto \int^t_0\lambda_5(t,s)\dif s\in L^1([0,T])$
such that for all $t,s\in[0,T]$, $\om\in\Om$ and $x,y\in{\mF^0_S}$
\ce
\|C_2(t,s,\om,x)-C_2(t,s,\om,y)\|^2_\mH&\leq& \lambda_5(t,s)\cdot\rho(\|x-y\|^2_{\mF^0_S}),\\
\|D_2(t,s,\om,x)-D_2(t,s,\om,y)\|^2_{L_2(\mU,\mH)}
&\leq& \lambda_5(t,s)\cdot\rho(\|x-y\|^2_{\mF^0_S}),
\de
where $\rho$ is same as in $\mathbf{(HF2)}$.

\item  There exist a positive progressively measurable process $\l_6$ and a positive real
function $\lambda_7$ satisfying
$\int^t_0\left[\l_7(t,s)+\mE\lambda_6(t,s)\right]\dif s\leq c_0\lambda^{2/q_1}_2(t)$,
and $0\leq \zeta\in L^1(\fA)$
such that for all $t,s\in[0,T]$, $\om\in\Om$ and $x\in{\mF^0_S}$
\ce
\|C_1(s,\om,x)\|^2_\mH+\|D_1(s,\om,x)\|^2_{L_2(\mU,\mH)}&\leq&
c_0\lambda^{2/q_1}_1(s)\cdot(\zeta(s,\om)+\|x\|^2_{\mF^0_S}),\\
\|C_2(t,s,\om,x)\|^2_\mH+\|D_2(t,s,\om,x)\|^2_{L_2(\mU,\mH)}&\leq&
\lambda_6(t,s,\om)+\l_7(t,s)\cdot\|x\|^2_{\mF^0_S},
\de
where $\lambda_1$ and $q_1$ are same as in $\mathbf{(HF1)}$.
\end{enumerate}

\bd\label{d4}
An $\mH$-valued continuous $\cF_t$-adapted process $X $ on $[-S,T]$ is called a solution of Eq.(\ref{Eq3}),
if
$$
\mE\left(\|X\|^2_{\mF^T_S}\right)+\|X_\cdot(0)\|^{q_1}_{\mK_{2,1}}
+\|X_\cdot(0)\|^{q_2}_{\mK_{2,2}}<+\infty,
$$
and (\ref{Eq3}) holds in $\mX^*$ for all $t\in[0,T]$ a.s..
\ed
\bt\label{th13}
Under $\mathbf{(HF1)}$-$\mathbf{(HF4)}$, there exists a unique solution to
Eq.(\ref{Eq3}) in the sense of Definition \ref{d4}.
\et
\begin{proof}
Let $X^1_t\equiv X_0(\cdot)$. One constructs the following iteration sequence $X^n_t$ for $n\in\mN$:
\ce
X^{n+1}_t(0)&=&X_0(0)+\int^t_0A(s,X^{n+1}_s(0))\dif s\\
&&+\int^t_0C_1(s,X^{n}_s)\dif s+\int^t_0\int^s_0C_2(s,r,X^{n}_r)\dif r\dif s\\
&&+\int^t_0 D_1(s,X^{n}_s)\dif W(s)+\int^t_0\int^s_0D_2(s,r,X^{n}_r)\dif W_r\dif s\\
&=:&X_0(0)+\int^t_0A^n(s,X^{n+1}_s(0))\dif s+\int^t_0B^n(s)\dif W(s),
\de
where
\ce
A^n(s,x)&:=&A^n_1(s,x)+A_2(s,x)\in\mY\subset\mX^*,\\
A^n_1(s,x)&:=&A_1(s,x)+G^n(s)\mX^*_1,\\
B^n(s)&:=&D_1(s,X^{n}_s)\in L_2(\mU,\mH)
\de
and
$$
G^n(s):=C_1(s,X^{n}_s)+\int^s_0C_2(s,r,X^{n}_r)\dif r+\int^s_0D_2(s,r,X^{n}_r)\dif W_r.
$$
First of all, we clearly have for any $x,y\in\mX$ and $s\in[0,T]$
$$
2[x-y, A^n(s,x)-A^n(s,y)]_\mX \leq\lambda_0(s)\cdot\|x-y\|^2_\mH.
$$
Secondly, by $\mathbf{(H3)}$ we have for any $x\in\mX$ and $s\in[0,T]$
\ce
&&2[x, A^n(s,x)]_\mX+\|B^n(s)\|^2_{L_2(\mU,\mH)}\\
&=&2[x,A(s,x)]_\mX +2\<x,G^n(s)\>_\mH
+\|D_1(s,X^{n}_s)\|^2_{L_2(\mU,\mH)}\\
&\leq&-\sum_{i=1,2}\Big(\lambda_i(s)\cdot\|x\|^{q_i}_{\mX_i}\Big)
+(\lambda_3(s)+1)\cdot\|x\|^2_\mH+\xi(s)\\
&&+\|G^n(s)\|^2_\mH+\|D_1(s,X^{n}_s)\|^2_{L_2(\mU,\mH)}.
\de
Moreover, by the embedding $\mH\subset\mX^*_1$ and $\mathbf{(H4)}$
we have for any $x\in\mX$ and $s\in[0,T]$
\ce
\|A^n_1(s,x)\|_{\mX^*_1}&\leq& \|A_1(s,x)\|_{\mX^*_1}+\|G^n(s)\|_{\mX^*_1}\\
&\leq&\eta_1(s)\cdot\lambda^{1/q_1}_1(s)+c_{A_1}\cdot\lambda_1(s,\om)\cdot\|x\|^{q_1-1}_{\mX_1}
+c_0\|G^n(s)\|_{\mH}.
\de
Hence, $(A^n,B^n)$ satisfies $\sH(\lambda_0,\lambda_1,\lambda_2,\lambda_3+1,\xi^n,\eta^n_1,\eta_2,q_1,q_2)$, where
\ce
\xi^n(s)&:=&\xi(s)+\|D_1(s,X^{n}_s)\|^2_{L_2(\mU,\mH)}+\|G^n(s)\|^2_\mH\\
\eta^n_1(s)&:=&\eta_1(s)+c_0\cdot\lambda^{-1/q_1}_1(s)\cdot\|G^n(s)\|_{\mH}.
\de

Noting that by $\mathbf{(HF4)}$
\ce
\mE\|C_1(s,X^{n}_s)\|^2_\mH+\mE\|D_1(s,X^{n}_s)\|^2_{L_2(\mU,\mH)}
&\leq& c_0\lambda^{2/q_1}_2(s)\left(\mE\zeta(s)+\mE\|X^n_s\|^2_{\mF^0_S}\right)\\
\mE\left\|\int^s_0C_2(s,r,X^{n}_r)\dif r\right\|^2_\mH&\leq&c_0\lambda^{2/q_1}_2(s)
\left(1+\mE\|X^n\|^2_{\mF^s_S}\right)\\
\mE\left\|\int^s_0D_2(s,r,X^{n}_r)\dif W_r\right\|^2_\mH&\leq&c_0\lambda^{2/q_1}_2(s)
\left(1+\mE\|X^n\|^2_{\mF^s_S}\right),
\de
we have by $q_1\geq 2$ and Young's inequality (\ref{Young})
\ce
\mE\left(\|G^n(s)\|^{\frac{q_1}{q_1-1}}_{\mH}\right)&\leq&\left(\mE\|G^n(s)\|^2_\mH\right)^{\frac{q_1}{2(q_1-1)}}\\
&\leq& c_0\cdot\lambda_1^{\frac{1}{q_1-1}}(s)\cdot\left(1+\mE\|X^n\|^2_{\mF^s_S}\right)^{\frac{q_1}{2(q_1-1)}}\\
&\leq& c_0\cdot\lambda_1^{\frac{1}{q_1-1}}(s)\cdot\left(1+\mE\|X^n\|^2_{\mF^s_S}\right).
\de
Therefore,
\ce
\int^t_0\mE\xi^n(s)\dif s&\leq& c_0+c_0\int^t_0(\lambda_1(s)+1)\cdot\mE\|X^n_s\|^2_{\mF^0_S}\dif s,\\
\int^t_0\mE\left(|\eta^n_1(s)|^{\frac{q_1}{q_1-1}}\right)\dif s&\leq& c_0+c_0\int^t_0
\mE\|X^n\|^2_{\mF^s_S}\dif s.
\de
Thus, by Theorem \ref{th2} we have
\be
\mE\left(\|X^{n+1}_\cdot(0)\|^2_{\mF^t_0}\right)
&\leq& c_0\left(\mE\|X_0(0)\|^2_\mH+\int^t_0\mE\left(\xi^n(s)+|\eta^n_1(s)|^{\frac{q_1}{q_1-1}}
+|\eta_2(s)|^{\frac{q_2}{q_2-1}}\right)\dif s\right)\no\\
&\leq& c_0+c_0\int^t_0(\lambda_1(s)+1)\mE\|X^n\|^2_{\mF^s_S}\dif s,\label{Op1}
\ee
where $c_0$ is independent of $n$.

By induction methods and Theorem \ref{th1}, $\{X^n_t, n\in\mN\}$ are thus well defined.
Moreover, by Theorem \ref{th2} we have
$$
\|X^n_\cdot(0)\|^{q_1}_{\mK_{2,1}}+\|X^n_\cdot(0)\|^{q_2}_{\mK_{2,2}}\leq c_n<+\infty.
$$

Setting
$$
g_n(t):=\sup_{k\leq n+1}\mE\left(\|X^k\|^2_{\mF^t_S}\right),
$$
we then have by (\ref{Op1})
\ce
g_n(t)&\leq& 2\sup_{k\leq n+1}\mE\left(\|X^k_\cdot(0)\|^2_{\mF^t_0}\right)+2\mE\|X_0\|^2_{\mF^0_S}\\
&\leq&c_0+c_0\int^t_0(\lambda_1(s)+1)g_n(s)ds.
\de
Applying Gronwall's inequality yields
\be
\sup_{n\in\mN}\mE\left(\sup_{t\in[0,T]}\|X^{n}_t\|^2_{\mF^0_S}\right)
=\sup_{n\in\mN}\mE\left(\|X^{n}\|^2_{\mF^T_S}\right)<+\infty.\label{es}
\ee

Next, set $Z^{n,m}_t:=X^n_t-X^m_t$. By It\^o's formula, ($\mathbf{H2}$) and ($\mathbf{HF2}$), we have
\ce
\|Z^{n+1,m+1}_t(0)\|^2_\mH&=&2\int^t_0[Z^{n+1,m+1}_s(0), A^n(s,X^{n+1}_s(0))-A^m(s,X^{m+1}_s(0))]_\mX\dif s\\
&&+2\int^t_0\<Z^{n+1,m+1}_s(0),(B^n(s)-B^m(s))\dif W(s)\>_\mH\\
&&+\int^t_0\|B^n(s)-B^m(s)\|^2_{L_2(\mU,\mH)}\dif s\\
&=&2\int^t_0[Z^{n+1,m+1}_s(0), A(s,X^{n+1}_s(0))-A(s,X^{m+1}_s(0))]_\mX\dif s\\
&&+\int^t_0\|D_1(s,X^{n}_s)-D_1(s,X^m_s)\|^2_{L_2(\mU,\mH)}\dif s+\sum^4_{i=1}I^i_t\\
&\leq&\int^t_0\lambda_0(s)\cdot\|Z^{n+1,m+1}_s(0)\|^2_\mH\dif s\\
&&+\int^t_0\lambda_3(s)\cdot\rho(\|Z^{n,m}_s\|^2_{\mF^0_S})\dif s
+\sum^4_{i=1}I^i_t,
\de
where
\ce
I^1_t&:=&2\int^t_0\<Z^{n+1,m+1}_s(0),C_1(s,X^{n}_s)-C_1(s,X^{m}_s)\>_\mH\dif s\\
I^2_t&:=&2\int^t_0\<Z^{n+1,m+1}_s(0),\int^s_0 C_2(s,r,X^{n}_r)-C_2(s,r,X^m_r)\dif r\>_\mH\dif s\\
I^3_t&:=&2\int^t_0\<Z^{n+1,m+1}_s(0),(D_1(s,X^{n}_s)-D_1(s,X^m_s))\dif W(s)\>_\mH\\
I^4_t&:=&2\int^t_0\<Z^{n+1,m+1}_s(0),\int^s_0(D_2(s,r,X^{n}_r)-D_2(s,r,X^m_r))\dif W_r\>_\mH\dif s.
\de
By Burkholder's inequality and Young's inequality (\ref{Young}), we have
\ce
\mE\left(\|I^3_\cdot\|_{\mF^t_0}\right)&\leq& c_0\mE\left(\int^t_0\|Z^{n+1,m+1}_s(0)\|^2_\mH\cdot
\|D_1(s,X^{n}_s)-D_1(s,X^m_s)\|^2_{L_2(\mU,\mH)}\dif s\right)^{1/2}\\
&\leq&\frac{1}{2}\mE\|Z^{n+1,m+1}_\cdot(0)\|^2_{\mF^t_0}
+c_0\int^t_0\lambda_3(s)\cdot\mE\rho(\|Z^{n,m}_s\|^2_{\mF^0_S})\dif s
\de
and
\ce
\mE\left(\|I^4_\cdot\|_{\mF^t_0}\right)&\leq&\int^t_0\mE\|Z^{n+1,m+1}_s(0)\|^2_\mH\dif s
+\int^t_0\int^s_0\lambda_5(s,r)\cdot\mE\rho(\|Z^{n,m}_r\|^2_{\mF^0_S})\dif r\dif s\\
&\leq&\int^t_0\mE\|Z^{n+1,m+1}_\cdot(0)\|^2_{\mF^s_0}\dif s
+\int^t_0\mE\rho(\|Z^{n,m}\|^2_{\mF^s_0})\left(\int^s_0\lambda_5(s,r)\dif r\right)\dif s,
\de
where we have used that $\rho$ is increasing and
$$
\sup_{r\in[0,s]}\|Z^{n,m}_r\|^2_{\mF^0_S}=\|Z^{n,m}\|^2_{\mF^s_0}.
$$
Similarly,
\ce
\mE\left(\|I^1_\cdot\|_{\mF^t_0}\right)
&\leq&\int^t_0\mE\|Z^{n+1,m+1}_\cdot(0)\|^2_{\mF^s_0}\dif s
+\int^t_0\lambda_3(s)\cdot\mE\rho(\|Z^{n,m}\|^2_{\mF^s_0})\dif s,\\
\mE\left(\|I^2_\cdot\|_{\mF^t_0}\right)
&\leq&\int^t_0\mE\|Z^{n+1,m+1}_\cdot(0)\|^2_{\mF^s_0}\dif s
+\int^t_0\mE\rho(\|Z^{n,m}\|^2_{\mF^s_0})\left(\int^s_0\lambda_5(s,r)\dif r\right)\dif s.
\de
Combining the above calculations, we obtain
\ce
\mE\|Z^{n+1,m+1}_\cdot(0)\|^2_{\mF^t_0}&\leq&2\int^t_0(\lambda_0(s)+3)
\cdot\mE\|Z^{n+1,m+1}_\cdot(0)\|^2_{\mF^s_0}\dif s\\
&&+c_0\int^t_0\lambda_8(s)\cdot\mE\left(\rho(\|Z^{n,m}\|^2_{\mF^s_0})\right)\dif s,
\de
where $\lambda_8(s):=\lambda_3(s)+\int^s_0\lambda_5(s,r)\dif r$.

By Gronwall's inequality and Jensen's inequality, we have
\be
\mE\|Z^{n+1,m+1}_\cdot(0)\|^2_{\mF^t_0}
\leq c_0\int^t_0\lambda_8(s)\cdot\rho(\mE(\|Z^{n,m}\|^2_{\mF^s_0}))\dif s.\label{es2}
\ee
Now setting
$$
g(t):=\limsup_{n,m\rightarrow\infty}\mE\|Z^{n+1,m+1}_\cdot(0)\|^2_{\mF^t_0}=
\limsup_{n,m\rightarrow\infty}\mE\|Z^{n+1,m+1}\|^2_{\mF^t_0},
$$
we then get by (\ref{es}) and Fatou's lemma
$$
g(t)\leq c_0\int^t_0\lambda_8(s)\cdot\rho(g(s))\dif s.
$$
Using Lemma \ref{le1} yields that
$$
g(t)=0.
$$
Therefore, there is an $\mH$-valued continuous adapted process $X$ such that
$$
\lim_{n\rightarrow\infty}\mE\|X^{n}-X\|^2_{\mF^T_0}=0.
$$

It remains to show that $X_t$ is a solution in the sense of Definition \ref{d4}.
Let $\tilde X(t)$ solve the following equation(Theorem \ref{th1})
\ce
\tilde X(t)&=&X_0+\int^t_0A(s,\tilde X(s))\dif s
+\int^t_0C_1(s,X_s)\dif s+\int^t_0\int^s_0C_2(s,r,X_r)\dif r\dif s\\
&&+\int^t_0 D_1(s,s,X_s)\dif W(s)+\int^t_0\int^s_0D_2(s,r,X_r)\dif W_r\dif s.
\de
As in estimating (\ref{es2}), we can prove that
$$
\mE\|X^{n+1}_\cdot(0)-\tilde X(\cdot)\|^2_{\mF^t_0}\leq c_0\int^t_0\lambda_8(s)\cdot
\rho(\mE\|X^{n}_\cdot-X_\cdot(0)\|^2_{\mF^s_0})\dif s.
$$
Taking limits and by the dominated convergence theorem, we get
$$
\lim_{n\rightarrow\infty}\mE\|X^{n+1}_\cdot(0)-\tilde X(\cdot)\|^2_{\mF^t_0}=0.
$$
So, for each $t\in[0,T]$
$$
X(t)=\tilde X(t), \quad a.s..
$$
Moreover, by Theorems \ref{th1} and \ref{th2}, for almost all $\om$,
$t\mapsto \tilde X(t,\om)$ is continuous in $\mH$, and
$$
\|\tilde X_\cdot(0)\|^{q_1}_{\mK_{2,1}}
+\|\tilde X_\cdot(0)\|^{q_2}_{\mK_{2,2}}<+\infty.
$$
The uniqueness follows from the similar calculations, and the proof is thus complete.
\end{proof}

We now consider the following stochastic Volterra evolution equation:
\be
X_t(0)=X_{0}(0)+\int^t_{0}A(s,X_s)\dif s
+\int^t_{0}C(t,s,X_s)\dif s+\int^t_{0} D(t,s,X_s)\dif W_s,\label{Eq4}
\ee
where $X_{0}$ is an $\cF_{0}$-measurable $\mF^0_S(\mH)$-valued random variable and $A=A_1+A_2$,
\ce
A_i: [0,T]\times\Om\times\mX_i&\to&\mX^*_i,\ \ i=1,2,\\
C:[0,T]\times[0,T]\times\Omega\times \mF^0_S(\mH)&\to&\mH,\\
D:[0,T]\times[0,T]\times\Omega\times\mF^0_S(\mH)&\to& L_2(\mU,\mH)
\de
are progressively measurable.

We make the following assumptions:
\begin{enumerate}[(\bf $\mathbf{HV}$1)]
\item $X_0\in L^2(\Omega,\cF_0,P;\mF^0_S(\mH))$ and $(A,0)$ satisfies
$\sH(\lambda_0,\lambda_1,\lambda_2,\lambda_3,\xi,\eta_1,\eta_2, q_1,q_2)$,
where $\l_i\in L^1([0,T])$, $i=0,1,2,3$ are non-random strictly  positive functions, $q_i\geq 2$,
and $0\leq\eta_i\in L^{\frac{q_i}{q_i-1}}(\fA)$, $i=1,2$, and $0\leq\xi\in L^1(\fA)$.

\item $C(\cdot,s,\om,x)$ and $D(\cdot,s,\om,x)$ are differentiable with respect to the first variable
for all $s,\om,x$, and there exist a positive real function $\lambda_5$
satisfying $t\mapsto \int^t_0\lambda_5(t,s)\dif s\in L^1([0,T])$
and an increasing concave function $\rho$ satisfying (\ref{rho}) such that for all
$t,s\in[0,T]$, $\om\in\Om$ and $x,y\in{\mF^0_S}$
\ce
\|\p_tC(t,s,\om,x)-\p_tC(t,s,\om,y)\|^2_\mH&\leq& \lambda_5(t,s)\cdot\rho(\|x-y\|^2_{\mF^0_S}),\\
\|\p_tD(t,s,\om,x)-\p_t D(t,s,\om,y)\|^2_{L_2(\mU,\mH)}
&\leq& \lambda_5(t,s)\cdot\rho(\|x-y\|^2_{\mF^0_S}).
\de

\item There exist a positive real function $\lambda_3\in L^1([0,T])$
such that for all $(s,\om)\in[0,T]\times\Om$ and $x,y\in{\mF^0_S}$
\ce
\|C(s,s,\om,x)-C(s,s,\om,y)\|^2_\mH&\leq& \lambda_3(s)\cdot\rho(\|x-y\|^2_{\mF^0_S}),\\
\|D(s,s,\om,x)-D(s,s,\om,y)\|^2_{L_2(\mU,\mH)}
&\leq& \lambda_3(s)\cdot\rho(\|x-y\|^2_{\mF^0_S}),
\de
where $\rho$ is same as in $\mathbf{(HV2)}$.

\item There exist a positive progressively measurable process $\l_6$ and a positive real
function $\lambda_7$ satisfying
$\int^t_0\left[\l_7(t,s)+\mE\lambda_6(t,s)\right]\dif s\leq c_0\lambda^{2/q_1}_2(t)$, and $0\leq \zeta\in L^1(\fA)$
such that for all $t,s\in[0,T]$, $\om\in\Om$ and $x\in{\mF^0_S}$
\ce
\|C(s,s,\om,x)\|^2_\mH+\|D(s,s,\om,x)\|^2_{L_2(\mU,\mH)}&\leq&
c_0\lambda^{2/q_1}_2(s)\cdot(\zeta(s,\om)+\|x\|^2_{\mF^0_S}),\\
\|\p_tC(t,s,\om,x)\|^2_\mH+\|\p_tD(t,s,\om,x)\|^2_{L_2(\mU,\mH)}
&\leq& \lambda_6(t,s,\om)+\l_7(t,s)\cdot\|x\|^2_{\mF^0_S},
\de
where $\lambda_1$ and $q_1$ are same as in $\mathbf{(HV1)}$.
\end{enumerate}

\bd\label{d5}
An $\mH$-valued continuous $\cF_t$-adapted process $X$ on $[-S,T]$ is called a solution of Eq.(\ref{Eq4}),
if
$$
\mE\left(\|X\|^2_{\mF^T_S}\right)+\|X_\cdot(0)\|^{q_1}_{\mK_{2,1}}
+\|X_\cdot(0)\|^{q_2}_{\mK_{2,2}}<+\infty,
$$
and (\ref{Eq4}) holds in $\mX^*$ for all $t\in[0,T]$ a.s..
\ed
\bt\label{th33}
Under $\mathbf{(HV1)}$-$\mathbf{(HV4)}$, there exists a unique solution to Eq.(\ref{Eq4})
in the sense of Definition \ref{d5}.
\et
\begin{proof}
Noting that by stochastic Fubini's theorem, we have
\ce
\int^t_0 D(t,s,X_s)\dif W_s&=&\int^t_0(D(t,s,X_s)-D(s,s,X_s))\dif W_s+\int^t_0 D(s,s,X_s)\dif W_s\\
&=&\int^t_0\int^t_s\p_rD(r,s,X_s)\dif r\dif W_s+\int^t_0 D(s,s,X_s)\dif W_s\\
&=&\int^t_0\int^s_0\p_sD(s,r,X_r)\dif W_r\dif s+\int^t_0 D(s,s,X_s)\dif W_s
\de
and
\ce
\int^t_0 C(t,s,X_s)\dif s&=&\int^t_0\int^s_0\p_sC(s,r,X_r)\dif r\dif s+\int^t_0 C(s,s,X_s)\dif s.
\de
Solving Eq.(\ref{Eq4}) is then equivalent to solving  the following equation
\ce
X_t&=&X_0+\int^t_0A(s,X_s(0))\dif s+\int^t_0 C(s,s,X_s)\dif s+\int^t_0\int^s_0\p_sC(s,r,X_r)\dif r\dif s\\
&&+\int^t_0 D(s,s,X_s)\dif W_s+\int^t_0\int^s_0\p_sD(s,r,X_r)\dif W_r\dif s.
\de
The result now follows from Theorem \ref{th13}.
\end{proof}


\section{Applications}

In this section, we discuss two applications, which in particular cover
the examples given in the introduction.

\subsection{Stochastic Porous Medium Equations}

Let $\cO$ be a bounded open subset of $\mR^d$.
For $q\geq 2$, let $W^{1,q}_0(\cO)$ and $W^{-1,\frac{q}{q-1}}(\cO)$ be the usual
Sobolev spaces(cf. \cite{Ze,Sh}).

For $q\geq 2$, set
$$
\mX_1=\mX_2=\mX:=L^q(\cO), \ \ \mH:=W^{-1,2}(\cO),\ \ \mX^*:=(L^{q}(\cO))^*.
$$
The inner product in $\mH$ is given by
$$
\<u,v\>_\mH:=\int_{\cO}(-\Delta)^{-1/2}u(x)\cdot(-\Delta)^{-1/2}v(x)~\dif x,\ \
u,v\in\mH.
$$
Note that $-\Delta$ establishes an isomorphism between $W^{1,2}_0(\cO)$
and $W^{-1,2}(\cO)$. We shall identify $W^{1,2}_0(\cO)$ with the dual space $\mH^*$ of $\mH$,
and hence $\mH^*=W^{1,2}_0(\cO)\subset L^{\frac{q}{q-1}}(\cO)$.
Thus, we have the evolution triple(cf. \cite{Sh, Roe})
$$
\mX\subset\mH\simeq\mH^*\subset \mX^*
$$
where $\simeq$ is understood via $-\Delta$. Moreover, for $u\in\mX$ and $v\in\mX^*$
$$
[u,v]_\mX=\int_\cO u(x)\cdot(-\Delta)^{-1}v(x)\dif x.
$$

Let $\varphi$ be a real measurable function on $[0,T]\times\Omega\times\mR$
satisfying the following assumptions:
\begin{enumerate}[{\bf (HP1)}]
\item For each $r\in\mR$, $(t,\omega)\mapsto\varphi(t,\omega,r)$ is a measurable adapted process.
\item For each $(t,\omega)\in[0,T]\times\Omega$, $r\mapsto\varphi(t,\omega,r)$ is continuous.
\item There exist $q\geq 2$ and positive functions $\xi,\eta, \lambda\in L^1(\fA)$,
where $\lambda(t,\omega)>0$ for $\dif t\times\dif P$ almost all $(t,\omega)\in[0,T]\times\Omega$,
such that for all
$(t,\omega,r)\in[0,T]\times\Omega\times\mR$
$$
r\cdot\varphi(t,\omega,r)\geq \lambda(t,\omega)\cdot |r|^q-\xi(t,\omega),
$$
and
$$
|\varphi(t,\omega,r)|\leq \lambda(t,\omega)\cdot|r|^{q-1}+
\eta^{\frac{q}{q-1}}(t,\omega)\cdot\lambda^{\frac{1}{q}}(t,\omega).
$$
\item For all $(t,\omega)\in[0,T]\times\Omega$ and $r,r'\in\mR$
$$
(r-r')\cdot(\varphi(t,\omega,r)-\varphi(t,\omega,r'))\geq 0.
$$
\end{enumerate}

Define the evolution operator $A$ as follows: for $u\in\mX=L^p(\cO)$
$$
A(t,\omega,u):=\Delta\varphi(t,\omega,u).
$$
Then $A(t,\omega,u)\in\mX^*$ and for $u,v\in\mX$
\be
[v,A(t,\omega,u)]_\mX=-\int_\cO v(x)\cdot\varphi(t,\omega,u(x))\dif x.\label{PP6}
\ee

Consider the following stochastic porous medium equation
with constant diffusion coefficient and Dirichlet boundary conditions(cf. \cite{Sh,Roe})
\be
\left\{
\begin{array}{ll}
\dif u(t)=\Delta(\varphi(t,\omega,u(t)))\dif t + B\dif W(t),\\
u(t,x)=0,\ \ \forall x\in\p\cO,\\
u(0,x)=u_0(x)\in W^{-1,2}(\cO),
\end{array}
\right.\label{Eq9}
 \ee
where $B\in L_2(\mU,\mH)$. Of course, $B$ can be some random cylindrical function
or linear function in $u$. For simplicity, we do not discuss this case(see next subsection).

We now check the above $A$ satisfies {\bf (H1)}-{\bf (H4)}.

For {\bf (H1)}, it is direct by {\bf (HP2)}, {\bf (HP3)} and the dominated convergence theorem.

For {\bf (H2)}, we have by (\ref{PP6}) and {\bf (HP4)}
\ce
&&[u-v,A(t,\omega,u)-A(t,\omega,v)]_\mX\\
&=&-\int_\cO (u(x)-v(x))\cdot
(\varphi(t,\omega,u(x))-\varphi(t,\omega,v(x))\dif x\leq 0.
\de

For {\bf (H3)} and {\bf (H4)}, we have by (\ref{PP6}) and {\bf (HP3)}
\ce
[u,A(t,\omega,u)]_\mX
&=&-\int_\cO u(x)\cdot \varphi(t,\omega,u(x))\dif x\\
&\leq&-\lambda(t,\omega)\cdot\|u\|^q_{\mX}+\xi(t,\omega)\cdot\mathrm{Vol}(\cO),
\de
and
\ce
\|A(t,\omega,u)\|_{\mX^*}&=&\left(\int_\cO|\varphi(t,\omega,u(x))|^{\frac{q}{q-1}}
\dif x\right)^{\frac{q-1}{q}}\\
&\leq& \lambda(t,\omega)\cdot\|u\|^{q-1}_\mX+\eta^{\frac{q}{q-1}}(t,\omega)
\cdot\lambda^{\frac{1}{q}}(t,\omega)
\cdot\mathrm{Vol}(\cO).
\de

Hence, $A$ satisfies {\bf (H1)}-{\bf (H4)}, and Theorem \ref{th1}
can be used to this situation. In particular, Eq.(\ref{Eq9})
contains Eq.(\ref{porous}) as a special case with $\varphi(t,\omega,r)=|w_t(\omega)|\cdot |r|^{p-2}r$
and $\lambda(t,\omega)=|w_t(\omega)|$.

\subsection{Stochastic Reaction Diffusion Equations}
As in the previous subsection, let $\cO$ be a bounded open subset
of $\mR^d$. We denote by $l^2$ the usual Hilbert space of square summable real number sequences.

We are given three measurable mappings
\ce
(a_1,\cdots, a_d)=:a:\!\!&&[0,T]\times\Omega\times\cO\times\mR\to\mR^d\\
b:\!\!&&[0,T]\times\Omega\times\cO\times\mR\to\mR\\
(\sigma_1,\cdots,\sigma_j,\cdots)=:\sigma:\!\!&&[0,T]\times\Omega\times\cO\times\mR\to l^2,
\de
which satisfy that
\begin{enumerate}[{\bf (HR1)}]
\item For each $(x,r)\in\cO\times\mR$, $(t,\omega)\mapsto a(t,\omega,x,r)$,
$b(t,\omega,x,r)$ and $\sigma(t,\omega,x,r)$
are measurable adapted processes.
\item For each $(t,\omega,x)\in[0,T]\times\Omega\times\cO$,
$r\mapsto a(t,\omega,x,r), b(t,\omega,x,r)$ are continuous.
\item For all $(t,\omega,x)\in[0,T]\times\Omega\times\cO$, $r,r'\in\mR$ and $j=1,\cdots,d$
$$
(r-r')\cdot(a_j(t,\omega,x,r)-a_j(t,\omega,t,r'))\geq 0.
$$

\item There exist $q_1\geq 2$ and positive functions $\xi_1,\eta_1, \lambda_1\in L^1(\fA)$,
where $\lambda_1(t,\omega)>0$ for $\dif t\times\dif P$ almost all $(t,\omega)\in[0,T]\times\Omega$,
such that for all $(t,\omega,r)\in[0,T]\times\Omega\times\mR$
and $j=1,\cdots,d$
$$
r\cdot a_j(t,\omega,x,r)\geq \lambda_1(t,\omega)\cdot |r|^{q_1}-\xi_1(t,\omega),
$$
and
$$
|a_j(t,\omega,x,r)|\leq \lambda_1(t,\omega)\cdot|r|^{q_1-1}
+\eta^{\frac{q_1}{q_1-1}}_1(t,\omega)\cdot\lambda^{\frac{1}{q_1}}_1(t,\omega).
$$
\item $b$ satisfies {\bf (HR3)} and {\bf (HR4)} with  different constant $q_2\geq 2$
and functions $\lambda_2,\xi_2,\eta_2$.

\item There exist positive functions $\lambda_0,\lambda_3,\xi_3\in L^1(\fA)$
such that for all $(t,\omega,x)\in[0,T]\times\Omega\times\cO$ and $r,r'\in\mR$
$$
\|\sigma(t,\omega,x,r)\|^2_{l^2}\leq\lambda_3(t,\omega)\cdot|r|^2+\xi_3(t,\omega)
$$
and
$$
\|\sigma(t,\omega,x,r)-\sigma(t,\omega,x,r')\|^2_{l^2}\leq
\lambda_0(t,\omega)\cdot|r-r'|^2,
$$
where for some $c_\sigma>0$,
$$
0\leq\lambda_0(t,\omega)<c_\sigma\cdot(\lambda_1(t,\omega)\wedge\lambda_2(t,\omega)),
$$
and (\ref{Con}) holds, $\lambda_1$ and $\lambda_2$ are from {\bf (HR4)} and {\bf (HR5)}.
\end{enumerate}

We consider the following stochastic reaction diffusion equation with Dirichlet boundary conditions:
\be
\left\{
\begin{array}{ll}
\dif u(t,x)=\Big[\sum^d_{i=1}\p_i a_i(t,\omega,x,\p_i u(t,x))-b(t,\omega,x,u(t,x))\Big]\dif t\\
\quad\quad\quad\quad\ \ +\sum_{j=1}^\infty \sigma_j(t,\omega,x,u(t,x))\dif W_j(t),\\
u(t,x)=0,\ \ \forall x\in\p\cO,\\
u(0,x)=u_0(x)\in L^2(\cO),
\end{array}
\right.\label{Eq8}
\ee
where $W_j(t)=\<W(t),\ell_j\>_{\mU}$ and $\{\ell_j,j\in\mN\}$
is an orthogonal basis of $\mU$.

Let
$$
\mX_1:=W^{1,q_1}_0(\cO), \ \
\mX_2:=L^{q_2}(\cO),\ \
\mH:=L^2(\cO)
$$
and
$$
\mX_1^*:=W^{-1,\frac{q_1}{q_1-1}}(\cO),\ \ \mX_2^*:=L^{\frac{q_2}{q_2-1}}(\cO).
$$
If we identify $\mH^*$ with $\mH$, then
$$
\mX_1\subset\mH\simeq\mH^*\subset\mX^*_1, \ \ \mX_2\subset\mH\simeq\mH^*\subset\mX^*_2
$$
are two evolution triples.

Now define for $u,v\in \mX_1$
\be
[v,A_1(t,\omega,u)]_{\mX_1}:=-\sum^d_{i=1}
\int_\cO a_i(t,\omega,x,\p_i u(x))\cdot \p_i v(x)\dif x\label{PP2}
\ee
and for $u,v\in\mX_2$
\be
[v,A_2(t,\omega,u)]_{\mX_2}:=-\int_\cO b(t,\omega,x,u(x))\cdot v(x)\dif x.\label{PP3}
\ee
Clearly, for each $(t,\omega)\in[0,T]\times\Omega$ and
$u\in\mX_1$, $[\cdot, A_1(t,\omega,u)]_{\mX_1}\in\mX_1^*$ and for each
$u\in\mX_2$, $[\cdot, A_2(t,\omega,u)]_{\mX_2}\in\mX_2^*$.
Thus,
$$
A_1(t,\omega,\cdot): \mX_1\to\mX_1^*,\ \ A_2(t,\omega,\cdot):\mX_2\to\mX_2^*.
$$
Moreover, we also define for $u\in\mH=L^2(\cO)$
$$
B(t,\omega,u):=\sum_{j=1}^\infty \sigma_j(t,\omega,\cdot,u(\cdot))\cdot\ell_j\in L_2(\mU,\mH).
$$

We now check the above $A$ and $B$ satisfy {\bf (H1)}-{\bf (H4)}.

For {\bf (H1)}, it is direct by {\bf (HR2)}, {\bf (HR4)}, {\bf (HR5)}
and the dominated convergence theorem.

For {\bf (H2)}, we have by {\bf (HR3)}, {\bf (HR5)} and {\bf (HR6)}
\ce
&&2[u-v,A(t,\omega,u)-A(t,\omega,v)]_\mX
+\|B(t,\omega,u)-B(t,\omega,v)\|^2_{L_2(\mU,\mH)}\\
&=&-2\sum_{i=1}^d\int_\cO \Big(a_i(t,\omega,x,\p_i u(x))-a_i(t,\omega,x,\p_i v(x))\Big)
\cdot\p_i (u(x)-v(x))\dif x\\
&&-2\int_\cO\Big(b(t,\omega,x,u(x))-b(t,\omega,x,v(x))\Big)\cdot (u(x)-v(x))\dif x\\
&&+\int_\cO\sum_{j=1}^\infty
|\sigma_j(t,\omega,x,u(x))-\sigma_j(t,\omega,x,v(x))|^2\dif x\\
&\leq&0+0+\lambda_0(t,\omega)\cdot\|u-v\|^2_{\mH}.
\de

For {\bf (H3)}, we have by {\bf (HR4)}-{\bf (HR6)}
\ce
&&2[u,A(t,\omega,u)]_\mX+\|B(t,\omega,u)\|^2_{L_2(\mU,\mH)}\\
&=&-2\sum_{i=1}^d\int_\cO  a_i(t,\omega,x,\p_i u(x))
\cdot\p_i u(x)\dif x\\
&&-2\int_\cO b(t,\omega,x,u(x))\cdot u(x)\dif x
+\int_\cO\|\sigma(t,\omega,x,u(x))\|^2_{l^2}\dif x\\
&\leq&-2d\lambda_1(t,\omega)\int_\cO \sum_i|\p_iu(x)|^{q_1}\dif x
+d\cdot\xi_1(t,\omega)\cdot\mathrm{Vol}(\cO)\\
&&-2\lambda_2(t,\omega)\int_\cO |u(x)|^{q_2}\dif x
+\xi_2(t,\omega)\cdot\mathrm{Vol}(\cO)\\
&&+\lambda_3(t,\omega)\int_\cO|u(x)|^2\dif x+\xi_3(t,\omega)\cdot\mathrm{Vol}(\cO)\\
&\leq&-c_0\sum_{i=1,2}\lambda_i(t,\omega)\cdot\|u\|^{q_i}_{\mX_i}
+\lambda_3(t,\omega)\cdot\|u|^{2}_{\mH}+\tilde\xi(t,\omega),
\de
where $c_0>0$ only depends on $q_1,d$, and
$\tilde\xi:=\mathrm{Vol}(\cO)\cdot(d\cdot\xi_1+\xi_2+\xi_3)$.

For {\bf (H4)}, we have by (\ref{PP2}) and {\bf (HR4)}
\ce
\|A(t,\omega,u)\|_{\mX^*_1}&\leq&c_{q_1}\cdot\sum_i\left(\int_\cO
|a_i(t,\omega,x,\p_i u(x))|^{\frac{q_1}{q_1-1}}\dif x\right)^{\frac{q_1-1}{q_1}}\\
&\leq& c_{q_1}\cdot\lambda_1(t,\omega)\cdot\|u\|^{q_1-1}_{\mX_1}
+\eta^{\frac{q_1}{q_1-1}}_1(t,\omega)\cdot\lambda^{\frac{1}{q_1}}_1(t,\omega)\cdot\mathrm{Vol}(\cO),
\de
and by (\ref{PP3}) and {\bf (HR5)}
\ce
\|A(t,\omega,u)\|_{\mX^*_2}&\leq& \left(\int_\cO
|b(t,\omega,x,u(x))|^{\frac{q_2}{q_2-1}}\dif x\right)^{\frac{q_2-1}{q_2}}\\
&\leq& \lambda_1(t,\omega)\cdot\|u\|^{q_2-1}_{\mX_2}
+\eta^{\frac{q_2}{q_2-1}}_2(t,\omega)\cdot\lambda^{\frac{1}{q_2}}_2(t,\omega)\cdot\mathrm{Vol}(\cO),
\de

Hence, $(A,B)$  satisfies {\bf (H1)}-{\bf (H4)}.
Thus, we may use Theorem \ref{th1} to Eq.(\ref{Eq8}).
In particular, Eq.(\ref{rea}) is a special case of Eq.(\ref{Eq8}).
In fact, we may take
\ce
a_i(t,\omega,x,r)&:=&|w_t(\omega)|\cdot r,\ \ i=1,\cdots, d,\\
b(t,\omega,x,r)&:=&|w_t(\omega)|\cdot |r|^{p-2}r,\\
\sigma_1(t,\omega,x,r)&:=&\sqrt{|w_t(\omega)|}\cdot r,\ \ \sigma_i=0,\ \ i=2,\cdots.
\de
Then $\lambda_0(t,\omega)=\lambda_1(t,\omega)=\lambda_2(t,\omega)=|w_t(\omega)|$
and $\eta_i=\xi_j=0$, $i=1,2, j=1,2,3$. Moreover, it is clear that
(\ref{Con}) holds in this case.

\vspace{5mm}

{\bf Acknowledgements:}

This work was done while the author was a fellow of Alexander-Humboldt Foundation.
He would  like to thank Professor Michael R\"ockner, his host in Bielefeld University,  for
providing him their Lecture Notes \cite{Roe} and a stimulating environment. He is also grateful to
Professor Jiagang Ren, Dr. Wei Liu and Dr. Huijie Qiao for their valuable discussions.
The referee's very useful suggestions are also acknowledged.

\end{document}